\documentclass[11pt,letterpaper]{amsart}
\usepackage{amssymb,amsfonts,amsthm,array,epsfig,graphics,graphicx,tabularx}
\usepackage{hyperref}
\usepackage[all]{hypcap}
\usepackage{subfigure}
\usepackage{enumerate}

\usepackage[usenames,dvipsnames]{color}

\newenvironment{proo}[1][Proof]{\noindent {\bf #1~: }}{\hfill$\Box$\medskip}

\newtheorem{theorem}{Theorem}[section]
\newtheorem{prop}[theorem]{{Proposition}}

\newtheorem{lemm}[theorem]{{Lemma}}

\newtheorem{coro}[theorem]{{Corollary}}

\newtheorem{clai}[theorem]{{Claim}}

\newtheorem{rema}[theorem]{Remark}

\theoremstyle{definition}
\newtheorem{defi}[theorem]{{Definition}}

\theoremstyle{remark}

\title{Lyapunov graphs of Non-singular Morse-Smale flows on \(S^1 \times S^2\)}
\author{Fangfang Chen}
\date{\today}

\begin{document}
\maketitle

\begin{abstract}

 Following B. Yu's work on Lyapunov graphs of non-singular Smale flows on $S^1\times S^2$, we characterize the Lyapunov graphs of non-singular Morse-Smale flows  on $S^1 \times S^2$ 
 by using  filtrating neighborhoods as the local data attached to the vertices.  More precisely, we   determine which oriented graphs with vertices   labeled by filtrating neighborhoods can be
realized as such Lyapunov graphs. 
%
 
\end{abstract}

\section{Introduction}
\subsection{Historic remarks and the aim of the paper}
A \emph{Morse-Smale flow}   is a smooth flow whose chain recurrent set consists of finitely many  hyperbolic closed orbits and fixed points, and  satisfies the transversality condition \cite{Sm}.  If a Morse-Smale flow has no fixed point, then we call it a \emph{non-singular Morse-Smale flow}, abbreviated as an \emph{NMS flow}.

The existence  of NMS flows is closely related to the topology of the underlying manifold.
 Asimov \cite{As1}  proved that a closed $n$-manifold  ($n\geq 4$) admits an NMS flow
if and only if its Euler  characteristic is zero. In dimension three,
Morgan \cite{Mo} proved that an irreducible closed orientable $3$-manifold   admits an NMS flow if and only if this manifold is a graph manifold. 
The periodic orbits of NMS flows on $3$-manifolds have also been studied from the viewpoint of indexed links. See Wada \cite{Wa} for $S^3$ and Chen--Yu \cite{CY} for the extension to graph manifolds.

For a closed orientable $3$-manifold $M$ that may not be  irreducible, Morgan \cite{Mo} obtained that there exists a number $k\geq0$ such that $M\sharp k(S^1 \times S^2)$ admits an NMS flow, based on the handle decomposition of $M$. In particular,  $S^1\times S^2$ admits  NMS flows according to the genus-$1$ Heegaard splitting of $S^1 \times S^2$. 
Thus $S^1\times S^2$, as the simplest reducible closed orientable
$3$-manifold, is a natural case to study separately.

To study NMS flows on $S^1\times S^2$, 
Lyapunov
graphs provide an effective way to
describe  the dynamical relations among the periodic orbits. They were  first introduced by Franks in his paper \cite{Fr} in the study of non-singular Smale flows (abbreviated as   \emph{NS flows}) on \(S^{3}\). 
 For  a Lyapunov function \(f: M\rightarrow \mathbb{R}\) associated to a smooth flow, a \emph{Lyapunov graph} is an oriented graph by identifying each connected component of \(f^{-1}(c)\) to a point for each \(c \in \mathbb{R}\). Moreover,
each edge is oriented by the flow direction. Importantly, a  flow may admit different Lyapunov graphs.
Using the Lyapunov graph,  de Rezende and her coauthors classified the
Smale flows on \(S^3\) (\cite{D1}), smooth flows on surfaces (\cite{D2}), and gradient-like flows on manifolds (\cite{D3}, \cite{D4}). 

In \cite{Yu2012},  Yu characterized the Lyapunov graphs of NS flows on \(S^1 \times S^2\), where each   vertex corresponding to a saddle basic set is labeled with the suspension of a subshift of finite type with a matrix $A$.
Although NMS flows form a special class of NS flows,  their Lyapunov graphs require finer information than that provided by
the matrix $A$.

Let $\phi_t$ be an NMS flow on $S^1 \times S^2$, and $L$ be a  Lyapunov graph of $\phi_t$. 
Note that the small neighborhood of each vertex in $L$ corresponds, via the Lyapunov function, to a special neighborhood of a periodic orbit of $\phi_t$, called a filtrating neighborhood.  
A \emph{filtrating neighborhood} $N(\gamma)$ of a periodic orbit $\gamma$ of $\phi_t$ is a compact connected neighborhood of $\gamma$, such that     $\phi_t$ is transverse to $\partial N(\gamma)$,  $\gamma$ is the maximal invariant set of $\phi_t$ in $N(\gamma)$, and orbit intersections with $N(\gamma)$ are connected. 
In the terminology of Morgan \cite{Mo}, these neighborhoods are called \emph{fattened round handles}.
Cordero, Martínez Alfaro and Vindel studied several important classes of fattened round handles for NMS flows on $S^1 \times S^2$ \cite{CAV}, including those with purely toroidal boundary and those whose non-toroidal boundary components consist of one $2$-sphere and one bitorus.



In this paper, we  characterize the Lyapunov graphs of NMS  flows on $S^1 \times S^2$ by using  filtrating neighborhoods as the local data attached to the vertices.  Namely, we   determine which oriented graphs with vertices   labeled by filtrating neighborhoods can be
realized as such Lyapunov graphs. 
This requires the full list of filtrating
neighborhoods of saddle periodic orbits, which we establish in
Section~\ref{s.classfy FN}.







 \subsection{Main results}
 To state the main results of this paper clearly, we first introduce some definitions.
 
 Let $\phi_t$ be an NMS flow on $S^1\times S^2$ with  a Lyapunov function  $f$, and  $L$ be the Lyapunov graph associated to $f$. 
A \emph{regular level surface} is a connected component of $f^{-1}(c)$, where $c$ is a regular value of $f$.
Note that any two points in each edge of $L$ correspond  to two homeomorphic     regular level surfaces.
We assign a \emph{weight} to each edge of $L$, which is the genus of the regular level surface associated to that edge.

Let $v$ be a vertex of $L$. 
Denote by  \(e_{v}^{+}\) and \(e_{v}^{-}\)  the numbers of incoming and  outgoing  edges connected to \(v\), respectively. If \(e_{v}^{-}\cdot e_{v}^{+} \neq 0\), then   \(v\) is called a \emph{saddle} vertex. If \(e_{v}^{+}=0\) (resp. \(e_{v}^{-} =0\)),  \(v\) is called a \emph{source} (resp. \emph{sink}) vertex.   
Now, we introduce the  abstract version of weighted Lyapunov
graphs.

\begin{defi}\label{def:AWLG0}
An \emph{abstract weighted Lyapunov graph}, abbreviated as \emph{AWLG}, is a finite, connected, oriented   graph \(L\) which satisfies the following conditions.
\begin{enumerate}

\item \(L\) possesses no oriented cycles.
 \item Each  source or sink vertex is connected by exactly one edge, and the weight of this edge is $1$.
\item Each saddle vertex $v$ is assigned a label from the list of
filtrating-neighborhood types in Proposition~\ref{prop:N-classification}
and Proposition~\ref{fil nghd twist}. The  numbers $e_v^+$, $e_v^-$, and the weights of all edges connected to
 $v$  agree with those prescribed by the Lyapunov graph of the
filtrating neighborhood indicated by the label.

\end{enumerate}
\end{defi} 

Here,
  Proposition~\ref{prop:N-classification} and Proposition~\ref{fil nghd twist} give a complete classification of filtrating neighborhoods of saddle periodic orbits of NMS flows on $S^1 \times S^2$. These classifications are denoted as (A.1)--(A.5), (B.1)--(B.8), (C.1)--(C.6), and (D.1)--(D.2).

Our first main result (Theorem \ref{thm:allseparable0}) characterizes the AWLGs that can be realized as Lyapunov graphs of NMS flows on $S^1 \times S^2$ for which every regular level surface is separable.
Here, a regular level surface is called \emph{separable} if   it separates $S^1\times S^2$.


\begin{theorem}\label{thm:allseparable0}
 Let $L$ be an  AWLG. $L$  is associated with an NMS flow $\phi_t$ on $S^1 \times S^2$ such that each regular level surface is separable if and only if
the following conditions hold:
 \begin{enumerate}
   \item the first Betti number $\beta_1(L)$ of $L$ is equal to $0$;
\item each saddle vertex is labeled with one of  (A.3), (B.2), (B.7), (B.8), (C.1), (C.3),  (C.4) and  (D.1);
\item the total number of vertices labeled with  (B.7), (B.8), or (C.3)
is at most $1$.

 \end{enumerate}

 \end{theorem}

Theorems \ref{thm:inseparable torus0} and \ref{thm:inseparable not torus0} characterize the AWLGs that can be realized as Lyapunov graphs of NMS flows on $S^1 \times S^2$ for which an inseparable regular level surface exists.



\begin{theorem}\label{thm:inseparable torus0}
 Let $L$ be an  AWLG. $L$  is associated with an NMS flow $\phi_t$ on $S^1 \times S^2$ such that there exists an inseparable regular level surface  homeomorphic to $T^2$
 if and only if the following conditions hold:
 \begin{enumerate}
   \item $\beta_1(L)=1$; 
\item each saddle vertex is labeled with one of  (A.3), (B.2),  (C.1),  (C.4) and  (D.1).

 \end{enumerate}

 \end{theorem}

\begin{theorem}\label{thm:inseparable not torus0}
 Let $L$ be an  AWLG. $L$  is associated with an NMS flow $\phi_t$ on $S^1 \times S^2$ 
 such that there exists an  inseparable  regular level surface  not homeomorphic to $T^2$
 if and only if the following conditions hold:
 \begin{enumerate}
   \item $\beta_1(L)=1$; 
   \item 
every edge in the cycle $C$ of $L$
has weight $0$ or $2$, and every edge not  in $C$ has
weight $1$;
\item each saddle vertex is labeled with one of the types listed in 
    Proposition \ref{prop:N-classification} and Proposition \ref{fil nghd twist},
    except for types (B.7), (B.8), and (C.3).

 \end{enumerate}

 \end{theorem}

\subsection{Organization of the paper}

In Section \ref{s.Pre}, we introduce some definitions  and elementary properties for filtrating neighborhoods and weighted Lyapunov graphs.
In Section \ref{s.classfy FN}, we classify the 
filtrating neighborhoods of saddle periodic orbits of  NMS flows on $S^1 \times S^2$. 
In Section \ref{s.LYGallseparable}, we discuss  which abstract weighted Lyapunov graphs can be realized as Lyapunov graphs of  NMS flows $\phi_t$ on $S^1 \times S^2$ for which every regular level surface is separable.
In Section \ref{s.LYG an inseparable}, we discuss  which abstract weighted Lyapunov graphs can be realized as Lyapunov graphs of  NMS flows $\phi_t$ on $S^1 \times S^2$ for which an inseparable regular level surface exists.


\section{Preliminaries}\label{s.Pre}

\subsection{Filtrating neighborhoods}

Let $\phi_t$ be an NMS flow on a   compact orientable $3$-manifold $M$, and $\gamma$ be a periodic orbit of $\phi_t$.

\begin{defi}
A \emph{filtrating neighborhood} $N(\gamma)$ of $\gamma$ is a connected neighborhood of $\gamma$ in $M$, such that:
\begin{enumerate}

\item $\phi_t $ is transverse to $\partial N(\gamma)$;
\item $\gamma$ is the maximal invariant set of $\phi_t$ in $N(\gamma)$;

\item
the intersection of any orbit of $\phi_t$ with $N(\gamma)$ is connected.

\end{enumerate}
 \end{defi}

\begin{rema}
 Filtrating neighborhoods are called    fattened round  handles by  Morgan  \cite{Mo}.

\end{rema}

 If $\gamma$ is a repelling (resp.  attracting) periodic orbit of $\phi_t$, then the filtrating neighborhood  $N(\gamma)$ of $\gamma$ is a solid torus, and $\phi_t|_{N(\gamma)}$ is topologically equivalent to the flow on $D^2 \times S^1 = \{(z,\theta)| z\in \mathbb C, |z|\leq 1,   \theta\in \mathbb R/2\pi\mathbb Z\}$  induced by the vector field $(\dot z,  \dot \theta ) = (z, 1)$ (resp. $ (\dot z,  \dot \theta ) = (-z, 1)$).

From now on, we discuss the flow on the filtrating neighborhood $N(\gamma)$ when $\gamma$ is a saddle periodic orbit of $\phi_t$.
Let $V$ and $\widetilde{V}$ be two solid tori with coordinates
$(z,\theta)\in D^2\times \mathbb R/2\pi\mathbb Z$.
We define two smooth flows $Y_t$ and $\widetilde{Y}_t$ on $V$ and $\widetilde{V}$ respectively, as follows.
\begin{itemize}
    \item let   $Y_t$ be the flow on $V$ induced by the vector field $(\dot{z}, \dot{\theta})=(\bar{z}, 1)$; 
    \item let   $\widetilde{Y}_t$ be the flow  on $\widetilde{V}$ induced by the vector field $(\dot z,\dot\theta)=(e^{i\theta}\bar z,1)$.
    \end{itemize}

\begin{figure}[htbp]
\centering
\subfigure[$c_1, c_2, c_3, c_4$ are the dividing curves.]{\includegraphics[width=0.35\textwidth]{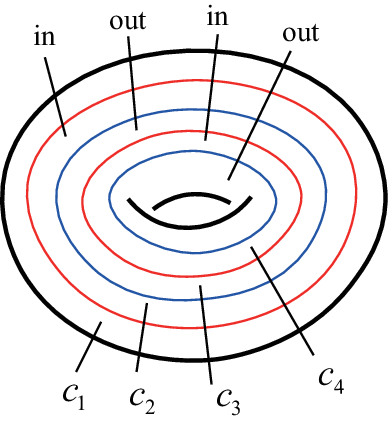}}
\hspace{.90in}
\subfigure[$c_1, c_2$ are the dividing curves.]{\includegraphics[width=0.35\textwidth]{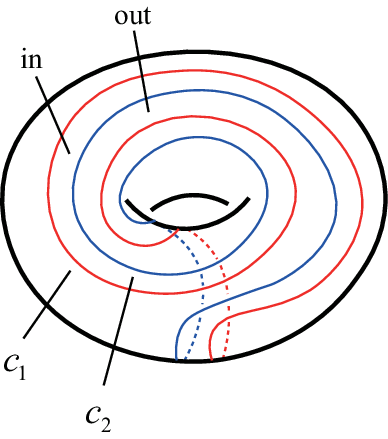}}

\caption{ }
\label{Fg_neig}
\end{figure}

     On $\partial V$,  the flow $Y_t$ is tangent to
$\partial V$ exactly along $4$ closed curves, such that these curves divide
$\partial V$ into $4$ annuli, and $Y_t$ is transverse to
the interiors of these annuli (see Figure \ref{Fg_neig} (a)). Similarly, on $\partial \widetilde{V}$, the flow  $\widetilde{Y}_t$ is tangent to
$\partial \widetilde{V}$ exactly along $2$ closed curves, such that these curves divide
$\partial \widetilde{V}$ into $2$ annuli, and  $\widetilde{Y}_t$ is transverse to
the interiors of these annuli (see Figure \ref{Fg_neig} (b)). We call these closed curves \emph{dividing curves}.
%

\par In \cite{Mo}, Morgan called \(V\) a \emph{kind orientable round 1-handle}, and called \(\widetilde{V}\) a \emph{kind non-orientable round 1-handle}.
In fact, if $\gamma$ is saddle, then there exists a tubular neighborhood $V(\gamma)$ of $\gamma$ such that $(V(\gamma), \phi_t|_{V(\gamma)})$ is topologically equivalent to either $(V, Y_t)$ or $(\widetilde{V}, \widetilde{Y}_t)$ \cite{PM}.
We call $\gamma$ a \emph{normal saddle periodic orbit} if $\gamma$ has a neighborhood which is topologically equivalent to $(V, Y_t)$. Otherwise, we call $\gamma$ a \emph{twisted saddle periodic orbit}.

\begin{figure}[htbp]
\centering
\includegraphics[scale=0.5]{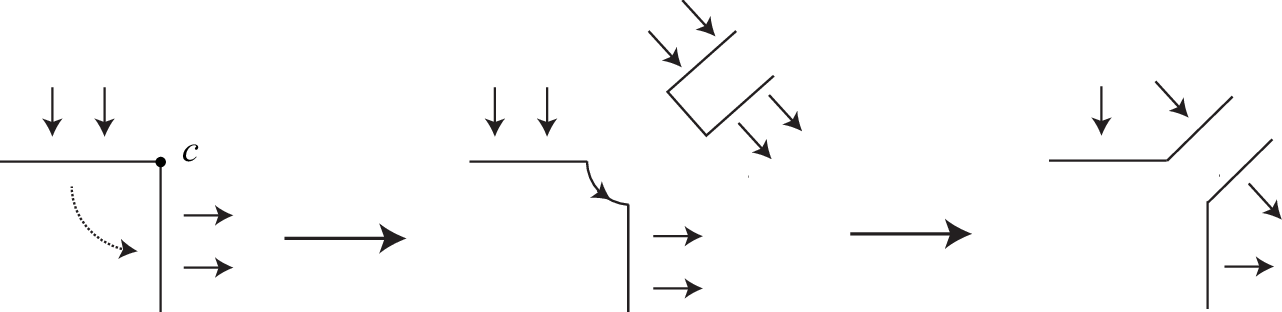}
\caption{Attaching a thickened surface to a dividing
curve $c$.}
\label{Fg_handle}
\end{figure}

Let \(\Sigma\) be a compact orientable surface with boundary. We call the $3$-manifold \(\Sigma \times I\) with vector field \(\frac{\partial}{\partial t}\) a \emph{thickened surface}, where $I=[0,1]$. 
Now, we  attach  a thickened surface $\Sigma \times I$  to $V$ or $\widetilde{V}$ along 
 a dividing
curve $c$  as follows.
 First, we blow up $c$ to obtain
   $c \times I$ and endow $c \times I$
 with the vector field $ \frac{\partial}{\partial t} $. Next, we choose a boundary component $c_0$ of $\Sigma$ and glue $c_0 \times I$ to $c\times I$ by preserving
the flowlines. Figure \ref{Fg_handle}  provides a sectional view of the attachment of a thickened surface. 
The following lemma comes from Theorem 4.4 in \cite{Yu}, which illustrates each filtrating neighborhood of a saddle periodic orbit can be obtained by attaching some thickened surfaces.

\begin{lemm}
\label{fil}
Suppose that $\gamma$ is a saddle periodic orbit of $\phi_t$, and let $N(\gamma)$ be a  filtrating neighborhood $N(\gamma)$ of $\gamma$. Then  $N(\gamma)$  always can be obtained by attaching some thickened surfaces along all dividing curves of $V$ or $\widetilde{V}$.
 \end{lemm}

\subsection{Weighted Lyapunov graphs}

The Lyapunov graph was first  used by  Franks  \cite{Fr} to classify nonsingular Smale flows on $S^3$. Let $\phi_t$ be an NMS  flow on  a closed orientable $3$-manifold $M$, and  let
$f:M\to \mathbb R$ be a Lyapunov function for $\phi_t$.   A \emph{Lyapunov graph} is an oriented  graph by identifying each connected component of \(f^{-1}(c)\) to a point for each \(c \in \mathbb{R}\),
where the components of the level sets of $f$  that contain periodic orbits produce the vertices of $L$.
Moreover,  each edge is oriented by the flow direction. 
A \emph{regular level surface} is a connected component of $f^{-1}(c)$, where $c$ is a regular value of $f$.
Note that  any two points in each edge of a Lyapunov graph correspond  to two homeomorphic     regular level surfaces.

Let $M'$ be a compact orientable $3$-manifold  with boundary, and $\psi_t$ be an NMS flow on $M'$ transverse to $\partial M'$. Suppose that $g$ is a Lyapunov function associated to $(M', \psi_t)$, and that $g$ is constant on each connected component of $\partial M'$. 
The Lyapunov graph of $(M',\psi_t)$ associated to $g$ is defined in the same way as     the case of $3$-manifolds without boundary. 
The difference is that the vertices consist of the vertices
corresponding to the periodic orbits of $\psi_t$, together with the endpoints corresponding to the connected
components of $\partial M'$.


\begin{defi}
A \emph{weighted Lyapunov graph} (abbreviated as \emph{WLG}) is a Lyapunov graph with a weight assigned to each edge, where the \emph{weight} of an edge is the genus of the regular level surface associated to that edge.
\end{defi}

\begin{prop}
\label{prop:regular level surface}
Let  $L$ be a WLG of  an NMS flow on $S^1 \times S^2$, then the weight of each edge in $L$ is $0,1,2$.
\end{prop}

Proposition \ref{prop:regular level surface} is due to Yu \cite{Yu2010}.The following lemma follows from the Poincaré--Hopf theorem \cite{Pu} for manifolds with boundary.

\begin{lemm}
\label{Euler}
Let \(\phi_{t}\) be a nonsingular smooth flow on a $3$-manifold  \(M\). Suppose that $\phi_t$ is  transverse to \(\partial M\), pointing inward along \(\partial^{+} M\) and outward along \(\partial^{-} M =\partial  M - \partial^{+} M\). Then \(\chi(\partial^{+} M) = \chi(\partial^{-} M)\), where \(\chi\) denotes the Euler characteristic.

\end{lemm}

Let $e$ be an edge of a WLG $L$. We denote by $w(e)$ the weight of $e$.
If $L - e$ is connected, we call $e$ \emph{inseparable}; otherwise, we call $e$ \emph{separable}.
Obviously, $e$ is inseparable if and only if the regular level surface associated to the midpoint of $e$ is inseparable.
Lemma \ref{lem. inseparable} is due to Proposition 4.4 of Yu \cite{Yu2012}.

\begin{lemm}
\label{lem. inseparable}
Let  \(L\) be a WLG of an NMS flow on \(S^{1} \times S^{2}\).
\begin{enumerate}
\item If each edge of $L$ is separable, then the weight of each edge of $L$ is $1$ and the first Betti number \(\beta_1(L)\) of $L$ is equal to $0$.

\item If at least one edge $e$ of $L$ is inseparable, then there are exactly the following two possibilities:
 \begin{enumerate}[(a)]
\item If   $w(e)=1$, then the weight of each edge in $L$ is $1$ and \(\beta_{1}(L)=1\).
\item Otherwise, $w(e)=0$ or $2$, and  \(\beta_{1}(L)=1\). Moreover, there exist at least one edge with  weight $0$, and  one edge with weight $2$.

\end{enumerate}
\end{enumerate}

\end{lemm}

\begin{prop}
\label{rek}
Let \(L\) be a  WLG of an NMS flow $\phi_t$ on \(S^{1} \times S^{2}\), and $e$ be an edge of \(L\) with weight   \(0\) or \(2\). Then
\begin{enumerate}
\item \(\beta_{1}(L)=1\);  
\item let \(C\) be a cycle in \(L\), then
every edge contained in $C$
has weight $0$ or $2$, and every edge not contained in $C$ has
weight $1$.

\end{enumerate}
\end{prop}

\begin{proo}
 By Lemma \ref{lem. inseparable}, we have \(\beta_{1}(L)=1\). Let \(C\) be a cycle in \(L\). If there exists an edge of \(C\) with weight \(1\), then by Lemma \ref{lem. inseparable},  the weight of each edge in $L$ is $1$, which contradicts that  \(w(e)\neq 1\). Thus the weight of each edge in \(C\) is \(0\) or \(2\).
 
 Let $e_0$ be an edge not contained  in $C$.
Cutting $L$ along the midpoint of $e_0$, we obtain two connected components $L_1$ and $L_2$, where $\beta_1(L_2)=0$. Let $M_2$ be the manifold corresponding to $L_2$ via the Lyapunov function, then $M_2$ has  only one boundary component and admits an NMS flow $\phi_t|_{M_2}$. By Lemma \ref{Euler}, $\partial M_2\cong T^2$, thus $w(e_0)=1$.

\end{proo}

Next, we discuss the embedding of the inseparable \(T^2\) and \(2T^2\) in \(S^1 \times S^2\), where $2T^2$ denotes a bitorus. Lemma  \ref{lem. T2} is due to Proposition 4.1 of Yu \cite{Yu2012}. 

\begin{lemm}
\label{lem. T2}
Let \(i: T^2  \hookrightarrow S^1 \times S^2\) be an embedding map, and let \(i_{\ast}: \pi_{1}(T^2) \to \pi_{1}(S^1 \times S^2)\) be the homomorphism induced by \(i\). Then:
\begin{enumerate}

\item If \(i(T^{2})\) is inseparable in \(S^{1}\times S^{2}\), then \(S^{1}\times S^{2} - i(T^{2}) \cong (O \sqcup K)^{c}\). Here \((O \sqcup K)^{c}\) is the complement of a two-component link in \(S^{3}\), where \(O\) is a trivial knot and \(K\) can be any knot unlinked with \(O\).

\item If \(i(T^{2})\) is separable in \(S^{1}\times S^{2}\), then \(i(T^2)\) bounds a knot complement in \(S^1 \times S^2 -i(T^2)\).
\end{enumerate}


\end{lemm}

\begin{prop}
\label{prop:2T2}
Let $\Sigma$ be an inseparable  $2T^2$ in $S^1 \times S^2$, then $S^1 \times S^2 - \Sigma$ is an open three manifold which can be embedded in $S^3$. 

\end{prop}

\begin{proo}
By Dehn's lemma \cite{Ha}, there is an essential simple closed curve \(c\) in $\Sigma$ which bounds a  disk $D$ in $S^1 \times S^2 - \Sigma$. Cutting $\Sigma$ along a small neighborhood of $\partial D$ and gluing two disks which is parallel to $D$, we obtain a new torus $T_0$ or two new tori $T_1$ and $T_2$, according as \(c\) is inseparable or separable in $\Sigma$.

 \begin{figure}[htbp]
\centering
\includegraphics[scale=0.5]{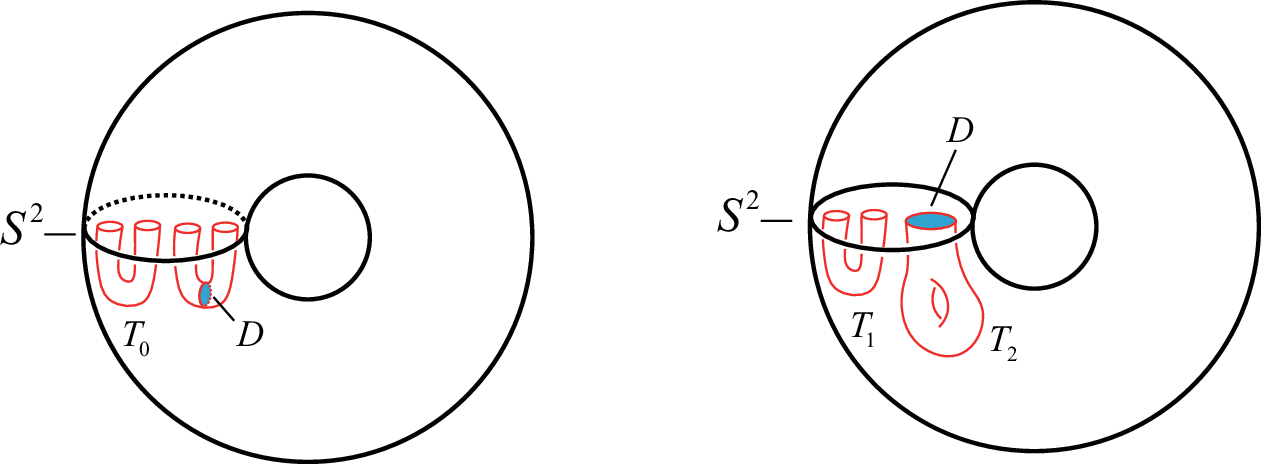}
\caption{ }
\label{Fg_2T2}
\end{figure}

 Then  $T_0$ is inseparable in \(S^{1} \times S^{2}\), or
one of \(T_{1}\) and \(T_{2}\) is inseparable in \(S^{1} \times S^{2}\). 
For the new inseparable torus, by applying Dehn's lemma again and gluing two disks, we can obtain an inseparable $2$-sphere in $S^1 \times S^2$.  Then the embedding of $T_0, T_1, T_2$ in $S^1 \times S^2$ is shown in Figure \ref{Fg_2T2}, where  we  assume that $T_1$ is inseparable for   the second case.
Then  
there exists a $3$-ball in $S^1 \times S^2 - T_1$ containing $T_2$. 
By Lemma   \ref{lem. T2} and the inverse of the surgery above, it is easy to observe that $S^1 \times S^2 - \Sigma$  can be embedded in $S^3$.

\end{proo}

\begin{figure}[htbp]
\centering
\includegraphics[scale=0.5]{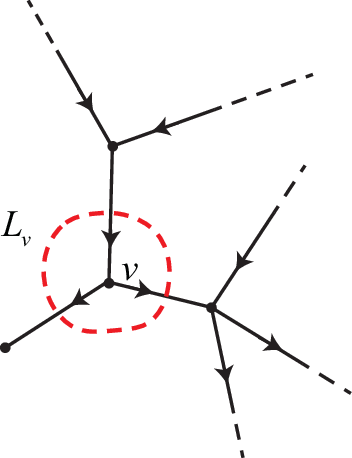}
\caption{$L_v$ is the star neighborhoods of the   vertex $v$.}
\label{Fg_Lyaneg}
\end{figure}

Let $L$ be a WLG of an NMS flow $\phi_t$ on a closed orientable $3$-manifold $M$.
  By cutting \(L\) along the midpoint of each edge connected to two vertices, we obtain some connected  graphs, each of which contains only one vertex of $L$. 
We call them \emph{the star neighborhoods of vertices} (see Figure \ref{Fg_Lyaneg}). 

Let $v_0$ be a vertex of $L$, corresponding to a periodic orbit
$\gamma_0$. The star neighborhood of $v_0$ corresponds, via the
Lyapunov function associated with $L$, to a filtrating neighborhood
$N_0$ of $\gamma_0$. Hence $L$ induces a decomposition of
$(M,\phi_t)$ into filtrating neighborhoods of periodic orbits. Such a
decomposition is called a fattened round handle decomposition by
Morgan~\cite{Mo}.

Conversely, let $\gamma_0,\gamma_1,\ldots,\gamma_n$ be all periodic
orbits of $\phi_t$, and let $N_0$ be a filtrating neighborhood of
$\gamma_0$. We can choose filtrating neighborhoods $N_i$ of
$\gamma_i$ for $i=1,\ldots,n$, such that
$N_0,N_1,\ldots,N_n$ are pairwise disjoint. For each
$j\in\{0,1,\ldots,n\}$, the entry and exit boundary components of
$\phi_t|_{N_j}$ determine a local weighted Lyapunov
graph. By gluing these local graphs according to the gluing rules of
filtrating neighborhoods, we obtain a WLG of $\phi_t$.

Hence, for NMS flows, the study
of WLGs can be reduced to the study of filtrating neighborhoods of periodic
orbits together with their gluing rules.


\begin{rema}\label{r.Ncondition}
Let $\phi_t$ be an NMS flow on $S^1 \times S^2$, and $N$ be a   filtrating neighborhood of a saddle periodic orbit of $\phi_t$.
By Proposition  \ref{prop:regular level surface} and Proposition \ref{rek},   $N$ should satisfy the following conditions:
\begin{enumerate}
    \item Each boundary component of $N$ is homeomorphic to $S^2$, $T^2$ or $2T^2$.
    \item If there exists a boundary component of $N$ that is not homeomorphic to $T^2$, then there are exactly two boundary components of $N$ not homeomorphic to $T^2$, and each of these two components is inseparable in $S^1\times S^2$.
    \item $N$ can be embedded in $S^1\times S^2$.
\end{enumerate}
\end{rema}
\section{Classification of filtrating  neighborhoods}\label{s.classfy FN}
\subsection{Normal  saddle periodic orbits}
\par Let $\gamma$ be a normal  saddle periodic orbit of an NMS flow $\phi_t$ on $S^1\times S^2$, and  $N$ be a filtrating  neighborhood of $\gamma$. Let $\partial^+N$ be the union of some boundary components of $N$ such that $\phi_t|_N$ enters $N$ transversely across $\partial^+N$, and let $\partial^-N = \partial N- \partial^+N$. Define
\[
W_1^s\sqcup  W_2^s = W^s(\gamma)\cap\partial^+ N,\qquad
W_1^u\sqcup  W_2^u = W^u(\gamma)\cap\partial^- N,
\]
where $W^s(\gamma)$ and $W^u(\gamma)$ denote the stable and unstable manifolds of $\gamma$ respectively.  

We say that $\gamma$ is \emph{trivial} in $N$  if $\gamma$ bounds a disk in $N$.
 $T^2_+,T^2_-$ denote tori, $S^2_+,S^2_-$ denote 2-spheres, and $2T^2_+,2T^2_-$ denote bitori. 

\begin{figure}[htbp]
\centering
\includegraphics[scale=0.67]{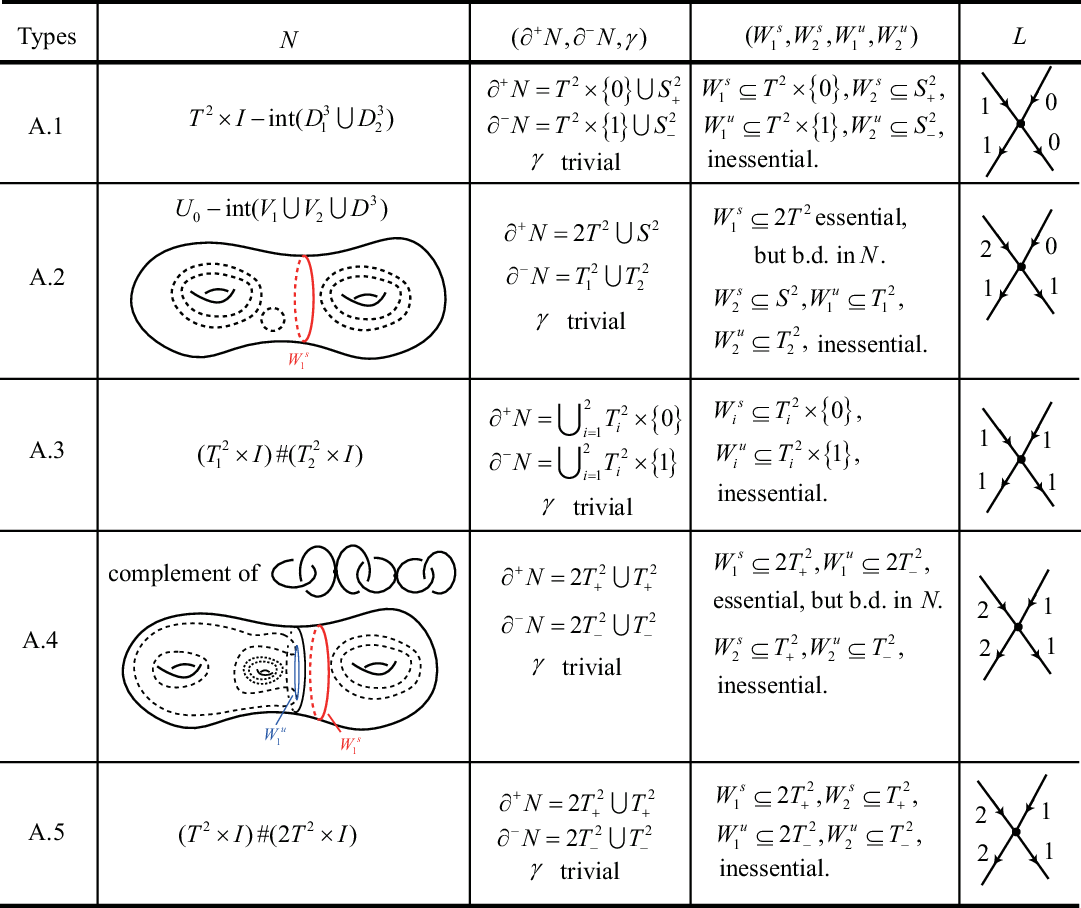}
\caption{Filtrating neighborhoods of types (A.1)--(A.5)}
\label{Fg_typeAnn}
\end{figure}

\begin{figure}[htbp]
\centering
\includegraphics[scale=0.67]{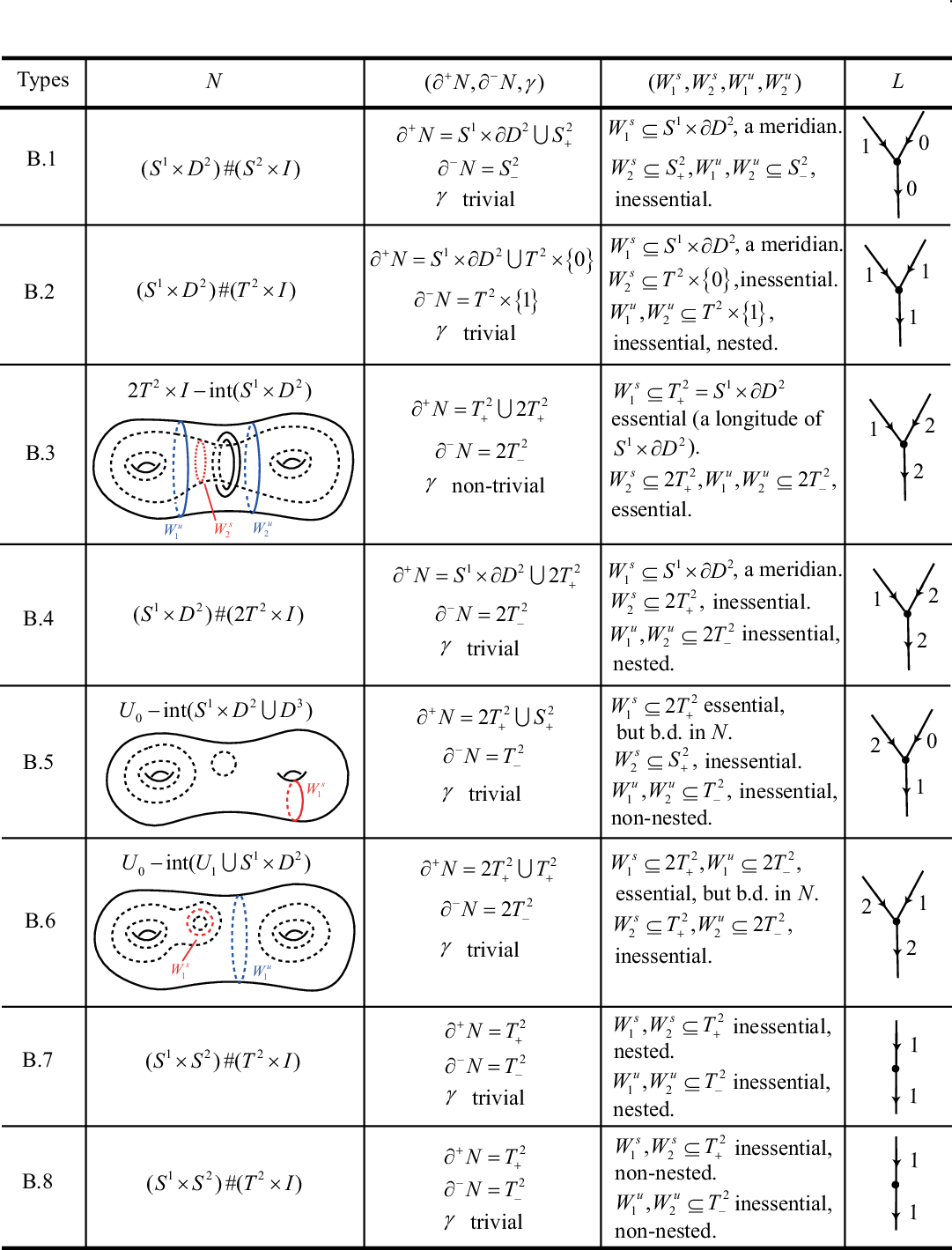}
\caption{Filtrating neighborhoods of types (B.1)--(B.8)}
\label{Fg_typeBnn}
\end{figure}

\begin{figure}[htbp]
\centering
\includegraphics[scale=0.67]{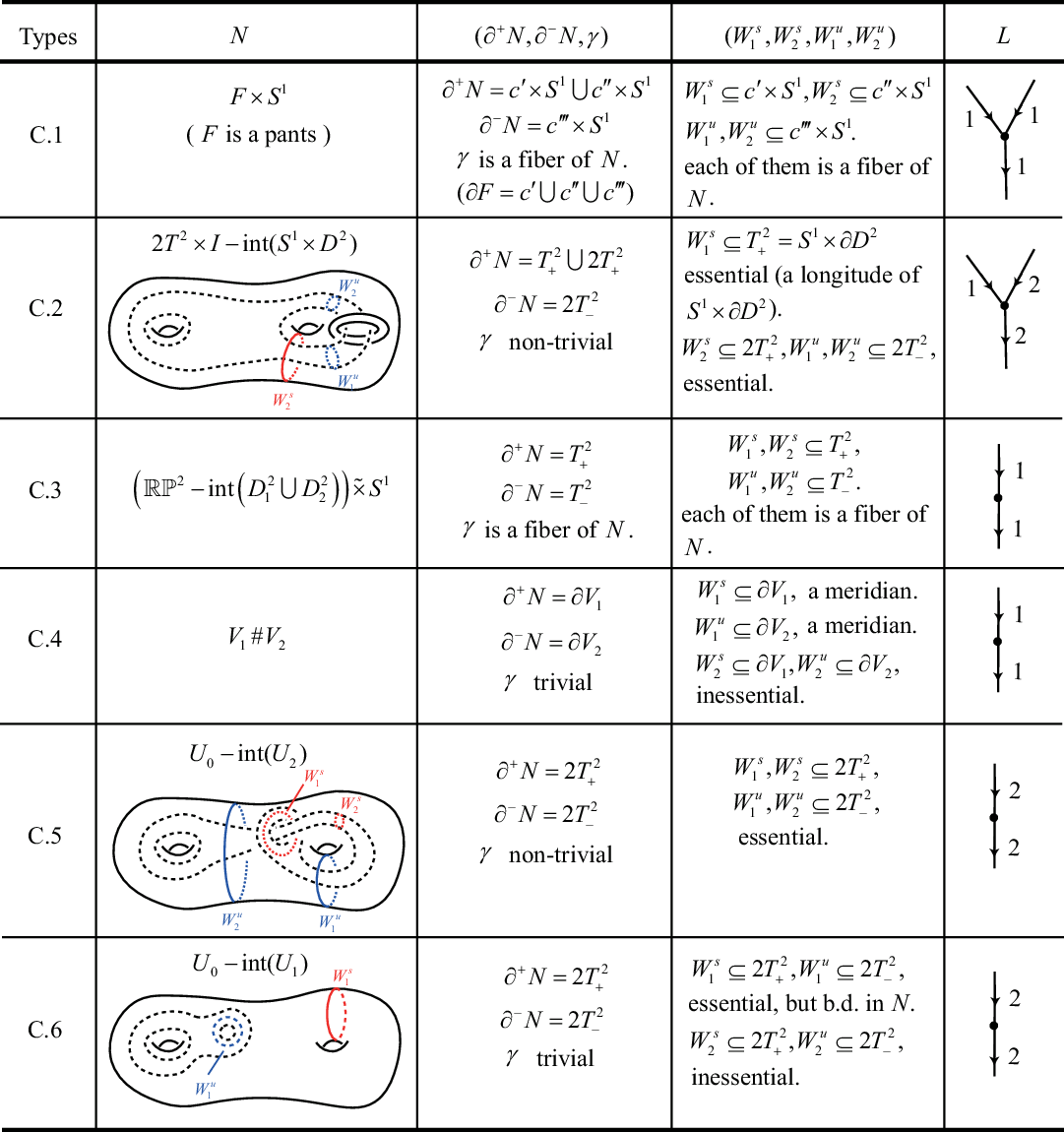}
\caption{Filtrating neighborhoods of types (C.1)--(C.6)}
\label{Fg_typeCnn}
\end{figure}

\begin{prop}\label{prop:N-classification}
Up to reversing time for the restricted flow $\phi_t|_N$, the filtrating neighborhood $N$ of the normal  saddle periodic orbit $\gamma$ is one of the following cases.
\begin{enumerate}[\indent \text{(A.1)}]

\item 
$N\cong T^2\times I - \operatorname{int}\bigl(D_1^3\cup  D_2^3\bigr)$
where $D_1^3,D_2^3$ are two disjoint  3-balls. $
\partial^+N = T^2\times\{0\}\cup S_+^2$, 
$\partial^-N = T^2\times\{1\}\cup S_-^2$, and 
 $\gamma$ is trivial in $N$. Moreover, $W_1^s\subseteq T^2\times\{0\}$ and $W_2^s\subseteq S_+^2$ are inessential. $W_1^u\subseteq T^2\times\{1\}$ and $W_2^u\subseteq S_-^2$ are inessential.
 \end{enumerate}
 
The remaining cases (A.2)--(A.5), (B.1)--(B.8) and (C.1)--(C.6) are illustrated in Figure  \ref{Fg_typeAnn}, Figure \ref{Fg_typeBnn} and Figure \ref{Fg_typeCnn}.
\end{prop}

\begin{rema}
\begin{enumerate}
\item
 In Figures  \ref{Fg_typeAnn}, \ref{Fg_typeBnn} and \ref{Fg_typeCnn}, 
 \emph{``b.d''} means that ``bounds a disk''.  we say that two inessential closed curves 
are
\emph{``nested''} if they bound two nested disks.
$U_0,U_1, U_2$ denote   genus two handlebodies, and $V_1, V_2$ denote  solid tori. Moreover,
$L$ denotes the weighted Lyapunov graph of $\phi_t|_{N}$.

\item For types (B.3), (C.2), (C.5), although $\gamma$ is not drawn in the figures, it has been reflected, because $\gamma$ is isotopic to $W_1^s, W_2^s, W_1^u, W_2^u$ in $N$.
\end{enumerate}
  
\end{rema}

\begin{proo}
By Remark  \ref{r.Ncondition}, $N$ should satisfy the following conditions:
\begin{enumerate}
    \item Each boundary component of $N$ is homeomorphic to $S^2$, $T^2$ or $2T^2$.
    \item If there exists a boundary component of $N$ that is not homeomorphic to $T^2$, then there are exactly two boundary components of $N$ not homeomorphic to $T^2$, and each of these two components is inseparable in $S^1\times S^2$.
    \item $N$ can be embedded in $S^1\times S^2$.
\end{enumerate}

By the condition (1), the attached thickened surface $\Sigma\times I$ satisfies
\[
\Sigma \cong S - \mathrm{int}\left(\bigcup_{i=1}^n D_i^2\right)
\]
where $S \cong S^2$, $T^2$ or $2T^2$, $D_i^2$ denotes a  disk, and $1\leq n\leq 3$.

Let $c_1, c_2, c_3, c_4$ be the dividing curves (see Figure \ref{Fg_neig} (a)). We use the notation $c_i+\Sigma$ to denote the operation of
attaching the thickened surface $\Sigma\times I$ along the dividing
curve $c_i$. For attachments along several dividing curves
$c_i,c_j,\ldots$, we write $c_i,c_j,\ldots+\Sigma$.

\textbf{Case 1. Attaching 4 thickened surfaces}
\begin{itemize}
    \item If $(c_1+D^2) \cup (c_2+D^2) \cup (c_3+D^2) \cup (c_4+D^2)$, then
  $\partial N$ consists of $4$ copies of $S^2$, which contradicts condition (2).
    
    \item $\bigl(c_1+\bigl(T^2-\mathrm{int}(D^2)\bigr)\bigr) \cup (c_2+D^2) \cup (c_3+D^2) \cup (c_4+D^2)$ corresponds to type (A.1).
    
    \item $\bigl(c_1+\bigl(T^2-\mathrm{int}\, D^2\bigr)\bigr) \cup \bigl(c_2+\bigl(T^2-\mathrm{int}\, D^2\bigr)\bigr) \cup (c_3+D^2) \cup (c_4+D^2)$ corresponds to type (A.2).
    
    \item $\bigl(c_1+\bigl(T^2-\mathrm{int}\, D^2\bigr)\bigr) \cup (c_2+D^2) \cup \bigl(c_3+\bigl(T^2-\mathrm{int}\, D^2\bigr)\bigr) \cup (c_4+D^2)$  corresponds to type (A.3).
    
    \item $\bigl(c_1 + \bigl(T^2 - \mathrm{int}\, D^2\bigr)\bigr) \cup \bigl(c_2 + \bigl(T^2 - \mathrm{int}\, D^2\bigr)\bigr) \cup \bigl(c_3 + \bigl(T^2 - \mathrm{int}\, D^2\bigr)\bigr) \cup (c_4 + D^2)$ corresponds to type (A.4).

    \item If $\bigcup_{i=1}^4 \bigl(c_i + \bigl(T^2 - \mathrm{int}\, D^2\bigr)\bigr)$, then $\partial N$ consists of four copies of $2T^2$, which contradicts condition (2).

    \item If $\bigl(c_1 + \bigl(2T^2 - \mathrm{int}\, D^2\bigr)\bigr) \cup (c_2 + D^2) \cup (c_3 + D^2) \cup (c_4 + D^2)$, then $N = 2T^2 \times I - \mathrm{int}\Bigl(\bigcup_{i=1}^2 D_i^2\Bigr)$, which contradicts condition (2).

    \item If $\bigl(c_1 + \bigl(2T^2 - \mathrm{int}\, D^2\bigr)\bigr) \cup \bigl(c_2 + \bigl(2T^2 - \mathrm{int}\, D^2\bigr)\bigr) \cup (c_3 + D^2) \cup (c_4 + D^2)$, then $N$ has a  boundary component homeomorphic to $4T^2$, which contradicts condition (1).

    \item If $\bigl(c_1 + \bigl(2T^2 - \mathrm{int}\, D^2\bigr)\bigr) \cup (c_2 + D^2) \cup \bigl(c_3 + \bigl(2T^2 - \mathrm{int}D^2\bigr)\bigr) + (c_4 + D^2)$, then $N = (2T^2 \times I) \# (2T^2 \times I)$, which contradicts condition (2).

    \item If $\bigl(c_1 + \bigl(T^2 - \mathrm{int}\, D^2\bigr)\bigr) \cup \bigl(c_2 + \bigl(2T^2 - \mathrm{int}\, D^2\bigr)\bigr) \cup (c_3 + D^2) \cup (c_4 + D^2)$, then $N$ has a  boundary component homeomorphic to $3T^2$, which contradicts condition (1).

    \item $\bigl(c_1 + \bigl(T^2 - \mathrm{int}\, D^2\bigr)\bigr) \cup (c_2 + D^2) \cup \bigl(c_3 + \bigl(2T^2 - \mathrm{int}\, D^2\bigr)\bigr) \cup (c_4 + D^2)$ corresponds to type (A.5).

    \item In all other cases, $N$ has  a boundary component homeomorphic to $mT^2$ ($m\ge 3$), which contradicts condition (1).
\end{itemize}

\textbf{Case 2. Attaching 3 thickened surfaces}

There exists one thickened surface satisfying
\[
\Sigma \cong S - \mathrm{int}\Bigl(\bigcup_{i=1}^2 D_i^2\Bigr), \quad \text{where } S = S^2,\ T^2 \text{ or } 2T^2.
\]
Since $N$ is a filtrating neighborhood, $\phi_t|_{N}$ is transverse to $\partial N$.
It is easy to observe that  $\partial \Sigma \times I$ is either normally glued to $N$ along two adjacent dividing curves, or non-normally glued to $N$ along two non-adjacent dividing curves.

\textbf{Case 2.1. Normal gluing}

\textbf{Case 2.1.1. $S \cong S^2$}
\begin{itemize}
    \item $(c_1,c_2 + \Sigma) \cup (c_3+D^2) \cup (c_4+D^2)$ corresponds to type (B.1).
    \item $(c_1,c_2 + \Sigma) \cup \bigl(c_3+(T^2-\mathrm{int}\,D^2)\bigr) \cup (c_4+D^2)$ corresponds to type (B.2).
    \item $(c_1,c_2 + \Sigma) \cup \bigl(c_3+(T^2-\mathrm{int}\,D^2)\bigr) \cup \bigl(c_4+(T^2-\mathrm{int}\,D^2)\bigr)$  corresponds to type (B.3).
    \item $(c_1,c_2 + \Sigma) \cup \bigl(c_3+(2T^2-\mathrm{int}\,D^2)\bigr) \cup (c_4+D^2)$  corresponds to type (B.4).
    \item In all other cases, $N$ has  a boundary component homeomorphic to $mT^2$ ($m\ge 3$), which contradicts condition (1).
\end{itemize}

\textbf{Case 2.1.2. $S \cong T^2$}

\begin{itemize}
    \item $(c_1,c_2 + \Sigma) \cup (c_3+D^2) \cup (c_4+D^2)$ corresponds to type (B.5).
    \item $(c_1,c_2 + \Sigma) \cup \bigl(c_3+(T^2-\mathrm{int}\,D^2)\bigr) \cup (c_4+D^2)$ corresponds to type (B.6).
    \item In all other cases, $N$ has  a boundary component homeomorphic to $mT^2$ ($m\ge 3$), which contradicts condition (1).
\end{itemize}

\textbf{Case 2.1.3. $S \cong 2T^2$}

It is easy to obtain that $N$ has  a boundary component homeomorphic to $mT^2$ ($m\ge 3$), which contradicts condition (1).

\textbf{Case 2.2. Non-normal gluing}

\textbf{Case 2.2.1. $S \cong S^2$}
\begin{itemize}
    \item If $(c_1,c_3 + \Sigma) \cup (c_2+D^2) \cup (c_4+D^2)$, then we can  observe that $(c_1,c_3+\Sigma)$ forms $S^1\times S^2 - \mathrm{int}N(c)$, where $N(c)$ is a tubular neighborhood of a closed braid (see Figure \ref{Fg_insepara} (a)). Thus $N \cong S^1\times S^2 - \mathrm{int}( D_1^3\cup D_2^3)$,
    and every component of $\partial N$ is separable in $S^1\times S^2$. This contradicts condition (2).

\begin{figure}[htbp]
\centering
\subfigure[ ]{\includegraphics[width=0.4\textwidth]{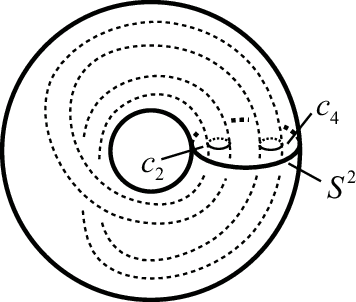}}
\hspace{.60in}
\subfigure[ ]{\includegraphics[width=0.34\textwidth]{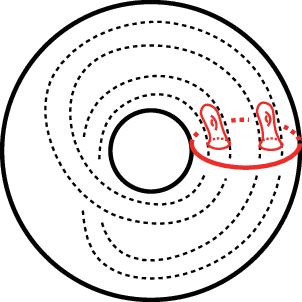}}

\caption{ }
\label{Fg_insepara}
\end{figure}

    \item $(c_1,c_3 + \Sigma) \cup \bigl(c_2+(T^2-\mathrm{int}\,D^2)\bigr) \cup (c_4+D^2)$ corresponds to type (B.7).

    \item If $(c_1,c_3 + \Sigma) \cup \bigl(c_2+(T^2-\mathrm{int}\,D^2)\bigr) \cup \bigl(c_4+(T^2-\mathrm{int}\,D^2)\bigr)$, then $\partial N$ consists of two copies of  $2T^2$. By condition (2), each boundary component is inseparable in $S^1\times S^2$. By Proposition \ref{prop:2T2}, $N$ can be embedded in $S^3$. However, $N$ contains an inseparable $2T^2$ (see Figure \ref{Fg_insepara} (b)), which contradicts the fact that   there is no inseparable closed surface in  $S^3$.

    \item If $(c_1,c_3 + \Sigma) \cup \bigl(c_2+(2T^2-\mathrm{int}\,D^2)\bigr) \cup (c_4+D^2)$, then by the same argument as the previous case, we obtain a contradiction.

    \item In all other cases, $N$ has  a boundary component homeomorphic to $mT^2$ ($m\ge 3$), which contradicts condition (1).
\end{itemize}

\textbf{Case 2.2.2. $S \cong T^2$}
\begin{itemize}
    \item $(c_1,c_3 + \Sigma) \cup (c_2+D^2) \cup (c_4+D^2)$ corresponds to type (B.8).

    \item If $(c_1,c_3 + \Sigma) \cup \bigl(c_2+(T^2-\mathrm{int}\,D^2)\bigr) \cup (c_4+D^2)$, then $\partial N$ consists of two copies of $2T^2$, and each component is inseparable in $S^1\times S^2$.  However,  $N$ contains an inseparable $T^2$, which contradicts  the conclusion from Proposition \ref{prop:2T2} 
  that $N$ can be embedded into $S^3$.

    \item In all other cases, $N$ has  a boundary component homeomorphic to $mT^2$ ($m\ge 3$), which contradicts condition (1).
\end{itemize}

\textbf{Case 2.2.3. $S \cong 2T^2$}
\begin{itemize}
    \item If $(c_1,c_3 + \Sigma) \cup (c_2+D^2) \cup (c_4+D^2)$, then $\partial N$ consists of two copies of $2T^2$, and $N$ contains an inseparable $2T^2$ (see Figure \ref{Fg_2T2insepa2}), which contradicts the embeddability of $N$ into $S^3$.

\begin{figure}[htbp]
\centering
\includegraphics[scale=0.5]{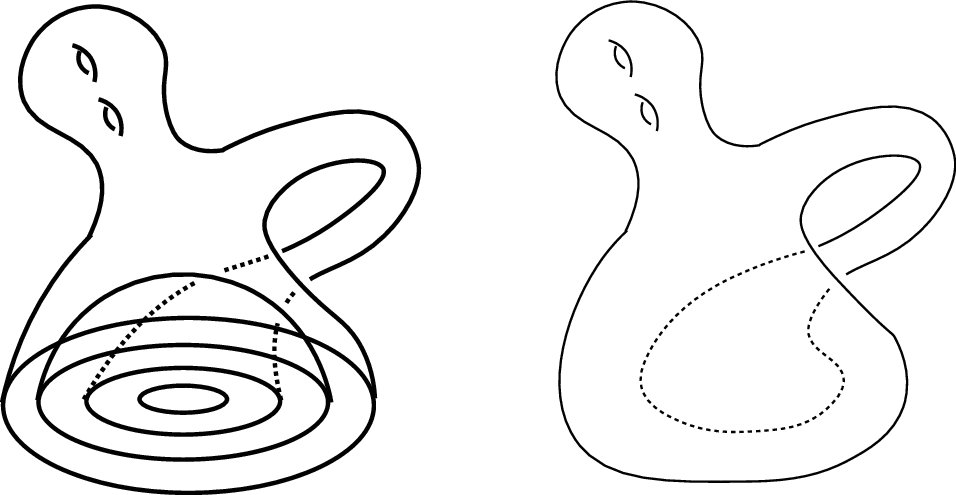}
\caption{On the left is $N$, and on the right is an inseparable $2T^2$ embedded in $N$.}
\label{Fg_2T2insepa2}
\end{figure}

    \item In all other cases, $N$ has  a boundary component homeomorphic to $mT^2$ ($m\ge 3$), which contradicts condition (1).
\end{itemize}

\textbf{Case 3. Attaching 2 thickened surfaces}

Let $\Sigma_j \cong S_j - \mathrm{int}\Bigl(\bigcup_{i=1}^{n_j} D_i^2\Bigr)$ be the two thickened surfaces ($j=1,2$),
where $n_1 = n_2 = 2$, or $n_1 = 3$ and $n_2 = 1$.

\textbf{Case 3.1. $n_1 = n_2 = 2$}

\textbf{Case 3.1.1. Normal gluing}
\begin{itemize}
    \item $\bigl(c_1,c_2 + (S^1\times I)\bigr) \cup \bigl(c_3,c_4 + (S^1\times  I)\bigr)$  corresponds to type (C.1).
    \item $\bigl(c_1,c_2 + (S^1\times I)\bigr) \cup \Bigl(c_3,c_4 + \bigl(T^2 - \mathrm{int}\bigl(\bigcup_{i=1}^2 D_i^2\bigr)\bigr)\Bigr)$  corresponds to type (C.2).
    \item In all other cases, $N$ has  a boundary component homeomorphic to $mT^2$ ($m\ge 3$), which contradicts condition (1).
\end{itemize}

\textbf{Case 3.1.2. Non-normal gluing}
\begin{itemize}
    \item $\bigl(c_1,c_3 + (S^1\times  I)\bigr) \cup \bigl(c_2,c_4 + (S^1\times I )\bigr)$  corresponds to type (C.3).

    \item If $\bigl(c_1,c_3 + (S^1\times  I)\bigr) \cup \Bigl(c_2,c_4 + \bigl(T^2 - \mathrm{int}\bigl(\bigcup_{i=1}^2 D_i^2\bigr)\bigr)\Bigr)$, then $\partial N$ consists of two copies of $2T^2$, and $N$ contains $T^2 \# 2\mathbb{RP}^2$, which contradicts the fact that $N$ can be embedded into $S^3$.

    \item In all other cases, $N$ has  a boundary component homeomorphic to $mT^2$ ($m\ge 3$), which contradicts condition (1).
\end{itemize}

\textbf{Case 3.2. $n_1 = 3$, $n_2 = 1$}
\begin{itemize}
    \item $\Bigl(c_1,c_2,c_3 + \bigl(S^2 - \mathrm{int}\bigl(\bigcup_{i=1}^3 D_i^2\bigr)\bigr)\Bigr) \cup (c_4 + D^2)$  corresponds to type (C.4).

    \item $\Bigl(c_1,c_2,c_3 + \bigl(S^2 - \mathrm{int}\bigl(\bigcup_{i=1}^3 D_i^2\bigr)\bigr)\Bigr) \cup \bigl(c_4 + (T^2 - \mathrm{int}\,D^2)\bigr)$  corresponds to type (C.5).

    \item $\Bigl(c_1,c_2,c_3 + \bigl(T^2 - \mathrm{int}\bigl(\bigcup_{i=1}^3 D_i^2\bigr)\bigr)\Bigr) \cup (c_4 + D^2)$  corresponds to type (C.6).

    \item In all other cases, $N$ has  a boundary component homeomorphic to $mT^2$ ($m\ge 3$), which contradicts condition (1).
\end{itemize}

By symmetry of the dividing curves, every admissible filtrating neighborhood is covered by the classification above, either for $\phi_t|_N$ or for its time reversal. Thus Proposition \ref{prop:N-classification} is proved.

\end{proo}

\begin{rema}

Up to topological equivalence, types (B.7) and (B.8) give two distinct filtrating neighborhoods.
For type (B.7), $W_1^s, W_2^s$ are nested on the entrance boundary component of the filtrating neighborhood, whereas this nesting property fails for type (B.8).

\end{rema}

\subsection{Twisted saddle periodic orbits}
\par Let $\gamma$ be a twisted saddle periodic orbit of an NMS flow $\phi_t$ on $S^1\times S^2$, and  $N$ be a filtrating  neighborhood of $\gamma$. Let 
$W^s = W^s(\gamma)\cap\partial^+ N$ and  $W^u= W^u(\gamma)\cap\partial^- N$. For knowledge of Seifert manifolds, please refer to Hatcher \cite{Ha}, and Martelli \cite{Ma}.

\begin{figure}[htbp]
\centering
\includegraphics[scale=0.5]{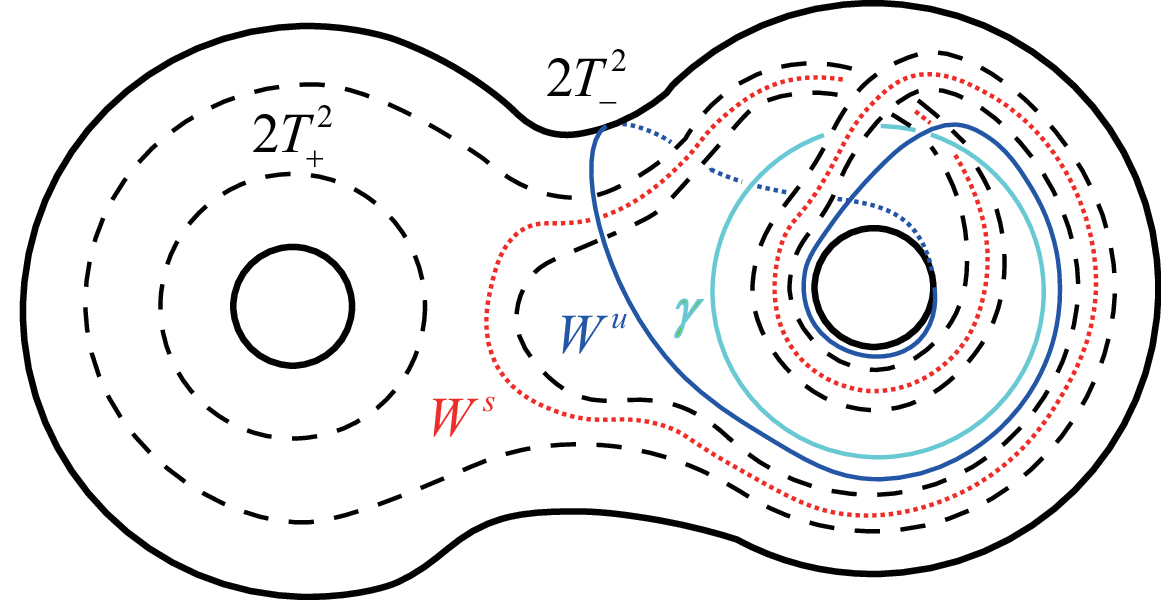}
\caption{The filtrating neighborhood of type (D.2)}
\label{Fg_D2}
\end{figure}

\begin{prop}
\label{fil nghd twist}
Up to reversing time for the restricted flow $\phi_t|_N$, the filtrating neighborhood $N$ of the twisted saddle periodic orbit $\gamma$ is one of the following cases.
\begin{enumerate}[\indent (D.1)]
  \item $N\cong  M(0,2;\frac{1}{2})$ is a Seifert manifold, whose base orbifold is an annulus with one cone points with the rotational isotropy group of order $2$.   $\partial^{+} N$ is a component of \(\partial N\). There is a Seifert fibering of $N$, such that $\gamma$ is  the singular fiber, and   $W^s,W^u$ are  two regular fibers.

  \item $N\cong U_0-\operatorname{int}U_3$, where   $U_0$ and $U_3$ are genus two handlebodies and $U_3$ is embedded in $U_0$. $\partial^{+} N$ is a component of \(\partial N\), and
  $\gamma, W^s,W^u$ are shown in Figure \ref{Fg_D2}.
\end{enumerate}
\end{prop}

\begin{proo}
Let $c_1$ and $c_2$ be the dividing curves of $\widetilde{V}$ (see Figure \ref{Fg_neig} (b)). Then there is either one   thickened surface or two thickened surfaces attached to $\widetilde{V}$. 

\textbf{Case 1. Attaching 2 thickened surfaces}

Let $\Sigma_i \times I$ be the thickened surface attached to $\widetilde{V}$ along the dividing curve $c_i$ ($i=1,2$).

\begin{itemize} 
    \item If $\Sigma_1 \cong D^2$ or $\Sigma_2 \cong D^2$, we may assume that $\Sigma_1 \cong D^2$. It is easy to observe that $c_1 + D^2$ yields $\mathbb{RP}^3 -  \mathrm{int}\,D^3$, which contradicts the embeddability of $N$ into $S^1\times S^2$.

    \item If $\Sigma_1 \cong \Sigma_2 \cong T^2 - \mathrm{int}\,D^2$, then $\partial N$ consists of two copies of $2T^2$, and $c_1+\Sigma_1$ yields $(T^2 \# \mathbb{RP}^2) \widetilde{\times} I$. This contradicts the conclusion from Proposition \ref{prop:2T2} that $N$ can be embedded in $S^3$.

    \item  In all other cases, $N$ has  a boundary component homeomorphic to $mT^2$ ($m\ge 3$), which is a contradiction.
\end{itemize}

\textbf{Case 2. Attaching one thickened surface}

Let $\Sigma \times I$ be the thickened surface attached to $\widetilde{V}$ along $c_1$ and $c_2$.
\begin{itemize}
    \item $\Sigma \cong S^1 \times I$  corresponds to type (D.1).
    \item $\Sigma \cong T^2 - \mathrm{int}(D_1^2\cup D_2^2)$  corresponds to type (D.2).
    \item  In all other cases, $N$ has  a boundary component homeomorphic to $mT^2$ ($m\ge 3$), which is a contradiction.
\end{itemize}

\end{proo}

\section{Lyapunov graphs with each regular level surface separable}\label{s.LYGallseparable}

Let $L$ be  an oriented  graph, and $v$ be a vertex of $L$. 
We denote by  \(e_{v}^{+}\) (resp. \(e_{v}^{-}\)) the number of incoming (resp. outgoing) edges connected to \(v\). If \(e_{v}^{-}\cdot e_{v}^{+} \neq 0\), then  we call \(v\) a \emph{saddle vertex}. If \(e_{v}^{+}=0\) (resp. \(e_{v}^{-} =0\)), we call \(v\) a \emph{source (resp. sink) vertex}. 

Suppose that  $L$ is  a WLG of an NMS flow $\phi_t$ on $S^1 \times S^2$.
Let $\gamma$ be the periodic orbit of $\phi_t$ corresponding to $v$, and  $N$ be a filtrating  neighborhood of $\gamma$ corresponding to   the star neighborhoods of $v$. Then the
following facts hold.
\begin{enumerate}

\item If   $v$ is a source or sink vertex, then $N \cong S^1 \times D^2$. Hence $v$ is connected by  exactly one edge,
and this edge has weight $1$.

\item  If  $v$ is a saddle   vertex, then $N$ is one of the types in Proposition \ref{prop:N-classification} and Proposition \ref{fil nghd twist}. We label the saddle vertex $v$ with the type of the  corresponding filtrating neighborhood $N$.
\end{enumerate}

  Now, we introduce the  abstract version of weighted Lyapunov
graphs.
\begin{defi}\label{def:AWLG}
An \emph{abstract weighted Lyapunov graph}, abbreviated as \emph{AWLG}, is a finite, connected, oriented   graph \(L\) which satisfies the following conditions.
\begin{enumerate}

\item \(L\) possesses no oriented cycles.
 \item Each  source or sink vertex is connected by exactly one edge, and the weight of this edge is $1$.
\item Each saddle vertex $v$ is assigned a label from the list of
filtrating-neighborhood types in Proposition~\ref{prop:N-classification}
and Proposition~\ref{fil nghd twist}. The  numbers $e_v^+$, $e_v^-$, and the weights of all edges connected to
 $v$  agree with those prescribed by the WLG of the
filtrating neighborhood indicated by the label.

 
\end{enumerate}
\end{defi}

\begin{rema}
The label of a saddle vertex $v$
corresponds to one of the local models listed in Proposition~\ref{prop:N-classification}
and Proposition~\ref{fil nghd twist}. 
Moreover, the star neighborhood of
$v$, together with the orientations and weights of its edges,
is isomorphic to the weighted Lyapunov graph of the corresponding local
model, possibly after reversing all
orientations. 

\end{rema}


\begin{defi}
Let $L$ be an AWLG. We say that $L$ \emph{is associated with} an
NMS flow on a $3$-manifold $M$ if there exists an NMS flow
$\phi_t$ on $M$ such that 
\begin{enumerate}
  \item$L$ is a WLG of $\phi_t$;
  \item for each saddle
vertex $v$, the filtrating neighborhood corresponding to the star neighborhood of $v$ has the same type as the label of $v$.
\end{enumerate}
\end{defi}

\begin{prop}
\label{S^3}
 Let $L$ be an  AWLG.  $L$ is associated with an NMS flow on $S^3$ if and only if the following conditions hold:
\begin{enumerate}
  \item $\beta_1(L)=0$;
  \item each saddle vertex is labeled with  one of (A.3), (B.2), (C.1), (C.4) and (D.1).
\end{enumerate}
\end{prop}

\begin{proo}
\emph{Necessity}.  By Franks \cite{Fr} and Wada \cite{Wa}, we obtain the necessity.

\par \emph{Sufficiency}. 
Let \(L\) be an AWLG satisfying the conditions (1) and (2), and
\(k\) be the number of saddle vertices of \(L\).
Obviously, \(L\) is associated with an NMS flow on \(S^3\) when \(k=0\).
 Assume  \(k \geq 1\),  and assume that an AWLG containing \(k-1\) saddle vertices, and  satisfying the conditions (1) and (2) is a WLG of an NMS flow on $S^3$.

We  cut \(L\) along an edge such that a component \(L_0\) of the resulting graph contains only one saddle vertex.  Let \(N\) be one of the  filtrating neighborhood of types (A.3), (B.2), (C.1), (C.4) and (D.1). Then we can always obtain a solid torus from \(N\) by attaching some solid tori suitably. In other words, \(L_0\)  can be realized as a WLG of an NMS flow on \(S^1 \times D^2\). By replacing \(L_0\) with a vertex, we obtain a new AWLG \(L'\) which has \(k-1\)  saddle vertices. By inductive assumption, \(L'\)  can be realized as a WLG of an NMS flow on $S^3$. Therefore, \(L\) is a WLG of an NMS flow on $S^3$.


\end{proo}

 Similar to the proof of the Sufficiency for Proposition \ref{S^3}, we can obtain the following corollary.
 
\begin{coro}\label{c.solid torus}
Let $L$ be an  AWLG, and $v_0$ be an endpoint of $L$. Suppose that  $\beta_1(L)=0$, and each saddle vertex is labeled with  one of (A.3), (B.2), (C.1), (C.4) and (D.1). Then 
$L$  is associated with an NMS flow  on $S^1 \times D^2$, such that $v_0$   corresponds to the boundary $S^1 \times \partial D^2$  via the Lyapunov
function.
\end{coro}

\begin{theorem}
 Let $L$ be an  AWLG. $L$  is associated with an NMS flow $\phi_t$ on $S^1 \times S^2$ such that each regular level surface is separable if and only if the following conditions hold:
 \begin{enumerate}
   \item $\beta_1(L)=0$ ;
\item each saddle vertex is labeled with one of  (A.3), (B.2), (B.7), (B.8), (C.1), (C.3),  (C.4) and  (D.1);
\item the total number of vertices labeled with  (B.7), (B.8), or (C.3)
is at most $1$.

 \end{enumerate}

 \end{theorem}

 \begin{proo}
\emph{Necessity}. Suppose that $L$ is a WLG of an NMS flow $\phi_t$ on $S^1 \times S^2$ such that each regular level surface is separable. By Lemma \ref{lem. inseparable}, every regular level surface of \((\phi_{t}, L)\)  is homeomorphic to \(T^{2}\) and \(\beta_1(L)=0\). By Proposition \ref{prop:N-classification} and Proposition \ref{fil nghd twist}, each saddle vertex of \(L\) is labeled with one of  (A.3), (B.2), (B.7), (B.8), (C.1), (C.3),  (C.4) and  (D.1).

\par Suppose that there exists a saddle vertex \(v\) labeled with   (B.7), (B.8) or (C.3) in $L$. Let $N$ be the  filtrating neighborhood corresponding to the star neighborhood of $v$. Then \(N \cong (S^1 \times S^2) \sharp (T^2 \times I)\)  or $N\cong M(-1,2;)$. Let  $T_1$ and $T_2$ be the boundary components of $N$. 
Cutting $S^1\times S^2$ along $T_1$,
we obtain two components $M_{1} \sqcup M_{2}$.
By Lemma \ref{lem. T2}, $M_{1}$ or $M_{2}$ is a knot complement. Without loss of generality, we assume that $M_{1}$ is a knot complement.
Since $N$  cannot be embedded in \(S^3\), we have
$N\subset M_2$.

Cutting $S^1\times S^2$ along $T_2$, we obtain two components $(M_1\cup_{T_1}N)\sqcup M_3.$
By Lemma \ref{lem. T2}, one of these two components is a knot
complement. Note that  $M_1\cup_{T_1}N$ cannot be embedded in $S^3$, then $M_3$ is a knot complement.
Therefore,
$S^1\times S^2-\operatorname{int}N=M_1\sqcup M_3$, where both $M_1$ and $M_3$ are knot complements. This implies that 
there is no more saddle vertices labeled with (B.7), (B.8) or (C.3) in $L$. 
\\
\par \emph{Sufficiency}. 
Let $L$ be an  abstract weighted Lyapunov graph satisfying conditions (1)-(3).
If there is no saddle vertex labeled with  (B.7), (B.8) or (C.3), then we can prove the sufficiency by induction, similar to the proof of Proposition \ref{S^3}. 

Suppose that there is a saddle vertex $v$ labeled with   (B.7), (B.8) or (C.3). then  $v$ is connected by exactly two edges.
Let $p_1$ and $ p_2$ be the midpoints of the two edges. Cutting $L$ along the points  $p_1$ and $p_2$, we obtain three components $L_0, L_1, L_2$, where $L_0$ is the  star neighborhood of $v$. Let $p_i^e$ be the endpoint of  $L_i$ related to $p_i$ for  $i=1,2$.
By Corollary \ref{c.solid torus}, each $L_i$  can be realized as a WLG of an NMS flow on a solid torus $V_i$, such that the endpoint $p_i^e$ corresponds to $\partial V_i$.

\par Suppose that \(v\) is labeled with (B.7) or  (B.8). Let \(N \cong (S^1 \times S^2) \sharp (T^2 \times I)\).
By gluing two solid tori $V_1, V_2$ to $N$ along \(\partial N\) suitably, we obtain \(S^1 \times S^2\).

\par Suppose that \(v\) is labeled with (C.3). Let $N\cong M(-1,2;)$. Namely, $N$ is the Seifert manifold whose base orbifold is  $\mathbb{RP}^2$ with two open disks removed.
Gluing the  solid torus $V_1$ to $N$ along one boundary component $T_1$ of \( N\), such that the Seifert fiber on $T_1$ bounds the meridian disk of $V_1$. Then we obtain \(V_{1} \cup_{T_1} N \cong (S^1 \times S^2 ) \sharp (S^1 \times D^2)\). 
Gluing $V_2$ to \(V_{1} \cup_{T_1} N\)  suitably,
 we obtain $S^1\times S^2$.

 Therefore, we can construct an NMS flow on $S^1 \times S^2$,  such that  $L$  is a weighted Lyapunov graph of this flow, and  each regular level surface is separable. The sufficiency is proved.

 \end{proo}

\section{Lyapunov graphs with  an inseparable regular level surface}\label{s.LYG an inseparable}
\begin{figure}[htbp]
\centering
\includegraphics[scale=0.5]{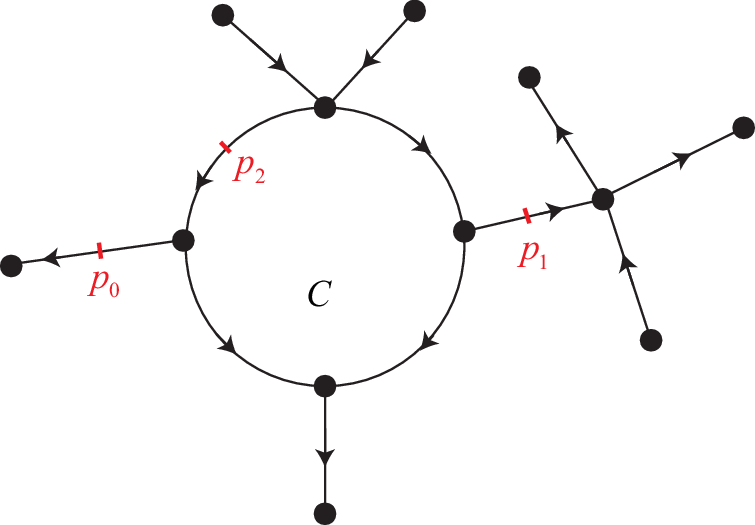}
\caption{$p_0 \notin \mathcal{P}_1(L)$, $p_1 \in \mathcal{P}_1(L)$, and $p_2 \in \mathcal{P}_2(L)$.}
\label{Fg_P1P2}
\end{figure}

Let $L$ be an  AWLG, such that any saddle vertex cannot be labeled with  (B.7), (B.8) or (C.3), and 
suppose that $\beta_1(L)=1$. Let $C$ be the cycle in $L$. Assume that each edge not in $C$ has weight  $1$.
Denote by $\mathcal{P}_1(L)$ the set of midpoints of edges of $L$
which have one endpoint in $C$ and the other endpoint at a saddle vertex
not lying in $C$.
Denote by $\mathcal{P}_2(L)$
the set of midpoints of the edges contained in $C$ (see Figure \ref{Fg_P1P2}).

Suppsose that $\mathcal{P}_1(L)=\varnothing$. Cutting $L$ along all points of $\mathcal{P}_2(L)$, we obtain some connected components, each of which contains exactly one saddle vertex. For each saddle vertex $v$, let $L_v$ be the component containing $v$. In fact,
$L_v$ is the union of the star neighborhood of $v$ and the star
neighborhoods of the sink and source vertices adjacent to $v$.

In the following, we define a compact 3-manifold $M(v)$ with boundary related to $v$.
We call $M(v)$ \emph{the saddle-vertex manifold associated with $v$}.

\begin{itemize}
    \item 
If $v$ is labeled with (A.3) or (B.2), and of the two edges in $C$ connected to $v$, one is incoming to $v$ and the other is outgoing from $v$, then  let $M(v)=T^2 \times I$.
\item
If $v$ is labeled with (A.3) or (B.2),  and  the two edges in $C$ connected to $v$ are both incoming to  $v$ or both outgoing from $v$,  then let $M(v)=(S^1 \times D^2) \sharp (S^1 \times D^2)$.
\item
If $v$ is labeled with (C.1) or (C.4),  then let $M(v)=(S^1 \times D^2) \sharp (S^1 \times D^2)$.

\item
If $v$ is labeled with (D.1), then let $M(v)=M(0,2;\frac{1}{2})$.
\end{itemize}

In the above cases, the weight of each edge connected to $v$ is $1$. Now, we suppose that there exists an edge connected to $v$ with a weight other than $1$.
Recall that each edge not in $C$ has weight $1$, then the 
  edge connected to $v$ without  weight $1$ must be contained in $C$.
In such cases, $M(v)$ is determined by the
weights of the two edges in $C$ connected to $v$ as follows.
\begin{itemize}

\item If the two edges in $C$ connected to $v$ both have weight
    $0$, then let
    $M(v)=S^2\times I$.

\item If the two edges in $C$ connected to $v$ both have weight
    $2$, then  
    $M(v)$ is of the form   $U-\operatorname{int} U_0$,
    where $U$ and $U_0$ are genus two handlebodies and $U_0$ is
    embedded in $U$.
    Moreover, $M(v)$ can be obtained from the standard local block $N$
corresponding to the label of $v$ by suitably gluing solid tori to all of the
torus boundary components of $N$.

 \item If one of the two edges in $C$ connected to $v$ has weight
    $0$, and the other has weight $2$, then let
    $M(v)=U-\operatorname{int} D^3$,
    where $U$ is a genus two handlebody and $D^3$ is a $3$-ball
    embedded in $U$.

%
%
\end{itemize}

\begin{rema}\label{r:S2,2T2boundary}
It follows from Proposition \ref{prop:N-classification} and Proposition \ref{fil nghd twist} that:

 \begin{enumerate}

\item if the two edges in $C$ connected to $v$ both have weight
    $0$ or both have weight $2$, then   one is incoming to $v$ and the other is outgoing from $v$;
    \item if one of the two edges in $C$ connected to $v$ has weight
    $0$, and the other has weight $2$, then  both edges are incoming to $v$ or both outgoing from $v$.

 \end{enumerate}

\end{rema}

\begin{lemm}\label{lem:block}
For each saddle vertex $v$,
there exists an NMS flow $\phi_t^v$ on $M(v)$, such that 
$L_v$ is a WLG of $\phi_t^v$, where the endpoints of $L_v$ related to the cutting points in $\mathcal{P}_2(L)$ correspond to the boundary of $ M(v)$ via the Lyapunov
function.
\end{lemm}

\begin{proo}
Suppose that \(v\) is labeled with (A.3). Then all possibilities of $L_v$ are shown in the Figure \ref{Fg_type(A.3)_LY}.
Let \(N \cong (T_1^2 \times I) \sharp (T_2^2 \times I)\) be the standard local
block of  type (A.3) in Proposition  \ref{prop:N-classification}, where $T_1, T_2$ are two tori. Moreover, $\partial^{+}N = (T_1^2 \times \{0\}) \sqcup (T_2^2 \times \{0\})$ and    $\partial^{-}N= (T_1^2 \times \{1\}) \sqcup (T_2^2 \times \{1\})$.

\begin{figure}[htbp]
\centering
\includegraphics[scale=0.8]{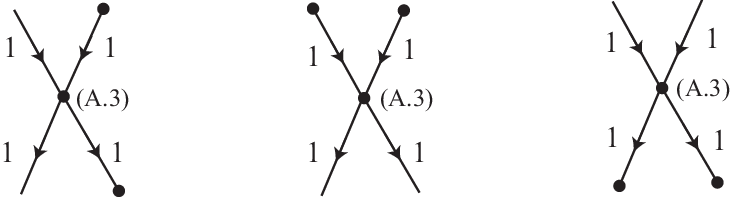}
\caption{$L_v$ for $v$  labeled with (A.3)}
\label{Fg_type(A.3)_LY}
\end{figure}

Depending on the orientation of the two edges of $L_v$ in $C$, we glue two solid tori to $N$  as follows, such that  the resulting manifold  is homeomorphic to $M(v)$.
\begin{itemize}

\item \textbf{One incoming, one outgoing:} we can obtain $T^2 \times I$  by gluing two solid tori to $N$ suitably  along $ T_1^2 \times \partial I$.

\item  \textbf{Both incoming or both outgoing:} we obtain $(S^1 \times D^2) \sharp (S^1 \times D^2)$ by gluing  two solid tori to $N$ along $(T_1^2 \times \{0\}) \sqcup (T_2^2 \times \{0\})$ or $(T_1^2 \times \{1\}) \sqcup (T_2^2 \times \{1\})$. 
\end{itemize}
 
 Based on the natural flow on $N$,  we can construct an NMS flow $\phi_t^v$ on $M(v)$, such that 
$L_v$ is a WLG of $\phi_t^v$, where the endpoints of $L_v$ related to the cutting points in $\mathcal{P}_2(L)$ correspond to the boundary of $ M(v)$ via the Lyapunov
function.
This proves the assertion for the case (A.3), and the other cases can be proved similarly.
\end{proo}

%
%
%
%

\begin{theorem}\label{thm:inseparable torus}
 Let $L$ be an  AWLG. $L$  is associated with an NMS flow $\phi_t$ on $S^1 \times S^2$ such that there exists an inseparable regular level surface  homeomorphic to $T^2$
 if and only if the following conditions hold:
 \begin{enumerate}
   \item $\beta_1(L)=1$; 
\item each saddle vertex is labeled with one of  (A.3), (B.2),  (C.1),  (C.4) and  (D.1).

 \end{enumerate}

 \end{theorem}

\begin{proo}
\emph{Necessity}.  Suppose that $L$ is a WLG of an NMS flow $\phi_t$ on $S^1 \times S^2$ such that  there exists an inseparable regular level surface $T$  homeomorphic to $T^2$. By Lemma \ref{lem. inseparable}, every regular level surface of \((\phi_{t}, L)\)  is homeomorphic to \(T^{2}\) and \(\beta_1(L)=1\). By Proposition \ref{prop:N-classification} and Proposition \ref{fil nghd twist}, each saddle vertex of \(L\) is labeled with one of  (A.3), (B.2), (B.7), (B.8), (C.1), (C.3),  (C.4) and  (D.1).

Let $v$ be a saddle vertex of $L$ labeled with  (B.7), (B.8), or  (C.3). Let   $N$ be the filtrating neighborhood corresponding to the star neighborhood of $v$. 
By Lemma \ref{lem. T2}, $S^1\times S^2 - T$ can be embedded in $S^3$, which implies that  $N$ can be embedded in $S^3$. 
However,  $N  \cong (S^1 \times S^2) \sharp (T^2 \times I)$ or $N\cong M(-1,2;)$, which is impossible. 
Thus  each saddle vertex is labeled with one of  (A.3), (B.2),  (C.1),  (C.4) and  (D.1).
\\
\par \emph{Sufficiency}. 
Let $L$ be an  AWLG satisfying conditions (1) and (2).
Suppose that $\mathcal{P}_1(L)\neq \varnothing$. Cutting $L$ along all points of $\mathcal{P}_1(L)$, we obtain some connected components $L_0, L_1, \dots, L_n$, where $\beta_1(L_0)=1$ and $\beta_1(L_i)=0$ for $i=1, \dots, n$. By Corollary \ref{c.solid torus}, $L_i$   can be realized as a WLG of an NMS flow on a solid torus, such that the end of $L_i$ related to the cutting point in  $\mathcal{P}_1(L)$ corresponds to the boundary of this solid torus.
By  replacing each $L_i$ with a vertex ($i=1, \dots, n$), we obtain a new AWLG $L'$ that still satisfies conditions (1) and (2), and with $\mathcal{P}_1(L')= \varnothing$.  Obviously, if $L'$   is associated with an NMS flow  on $S^1 \times S^2$, then so is $L$. Therefore,
from now on, we may assume without loss of generality that $\mathcal{P}_1(L)= \varnothing$.

Cutting $L$ along all points of $\mathcal{P}_2(L)$, we obtain some connected components, each of which contains exactly one saddle vertex. For each saddle vertex $v$, let $L_v$ be the component containing $v$.  
 As previously defined,  the saddle-vertex manifold
$M(v)$ associated with $v$ is homeomorphic to one of $T^2 \times I$, $(S^1 \times D^2) \sharp (S^1 \times D^2)$, and  $M(0,2;\frac{1}{2})$. Let $C$ be the unique cycle of $L$.

\begin{clai}\label{c:saddle revise orientation}
There exist two saddle vertices $v_1$ and  $v_2$ of $L$, such that 
\begin{enumerate}
   \item the two edges in $C$ connected to $v_1$ are both outgoing from $v_1$;
\item the two edges in $C$ connected to $v_2$ are both incoming to $v_2$;
\item $M(v_1)\cong M(v_2)\cong (S^1 \times D^2) \sharp (S^1 \times D^2)$.

\end{enumerate}

\end{clai}

\begin{proo}
According to the definition of $M(v)$,   if $M(v)$ is homeomorphic to $T^2 \times I$  or $M(0,2;\frac{1}{2})$, then one of the two edges in $C$ connected to $v$ 
is incoming to $v$ and the other is outgoing from $v$.
Then  if the two
edges in $C$ connected to $v$ are both incoming or both outgoing, then
$M(v)\cong (S^1 \times D^2)\sharp(S^1 \times D^2)$.

Since $\beta_1(L)=1$ and $L$ contains no oriented cycle,   there must exist at least one  saddle vertex $v_1$ such that the two edges in $C$ connected to $v_1$ are both outgoing from $v_1$, and at least one  saddle vertex $v_2$ such that the two edges in $C$ connected to $v_2$ are both incoming to $v_2$. Then $M(v_1)\cong M(v_2)\cong (S^1 \times D^2) \sharp (S^1 \times D^2)$, and  Claim \ref{c:saddle revise orientation} is proved.

\end{proo}

Let $m$ be the number of the saddle vertices of $L$. By Claim \ref{c:saddle revise orientation}, we have  $m\geq 2$. Let $v_3, \dots, v_m$ be other saddle vertices of $L$.
By Lemma \ref{lem:block}, for each $j=1, \dots, m$, there exists an NMS flow $\phi_t^{v_j}$ on $M(v_j)$, such that 
$L_{v_j}$ is a WLG of $\phi_t^{v_j}$, where the ends of $L_{v_j}$ related to the cutting points in $\mathcal{P}_2(L)$ correspond to the boundary of $ M(v_j)$ via the Lyapunov
function.

 \begin{clai}\label{c:inseparable torus gluing}
 By gluing the saddle-vertex manifolds $\{M(v_j)|j=1, \dots, m\}$ according to the gluing rules of $\{L_{v_j}| j=1,\dots, m\}$ in $L$, we can obtain $S^1\times S^2$.
 \end{clai}
 
 \begin{proo}
 If $m=2$, then by Claim \ref{c:saddle revise orientation}, we have $M(v_1)\cong M(v_2)\cong (S^1 \times D^2) \sharp (S^1 \times D^2)$. By gluing these two saddle-vertex manifolds along their boundary suitably, we can obtain \(S^1 \times S^2\).
 
Suppose that  $m>2$. Since $\beta_1(L)=1$, there is a saddle vertex in $\{v_3,\dots, v_m\}$  adjacent to $v_2$. We may assume that  $v_3$ is adjacent to $v_2$. Let $p_{(2,3)}$ be the point in $ \mathcal{P}_2(L)$ adjacent to both $v_2$ and $v_3$ (see Figure \ref{Fg_cut}).
\begin{figure}[htbp]
\centering
\includegraphics[scale=0.7]{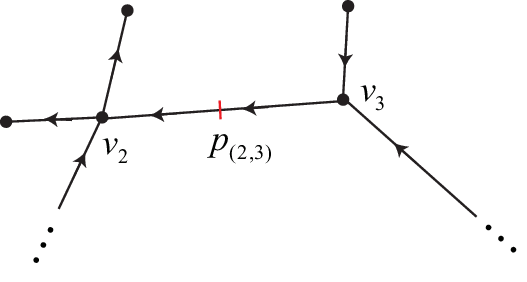}
\caption{ }
\label{Fg_cut}
\end{figure}

Then $M(v_3)$ is homeomorphic to one of $T^2 \times I$, $(S^1 \times D^2) \sharp (S^1 \times D^2)$, and $M(0,2;\frac{1}{2})$. 
We   glue $M(v_2)$ and $M(v_3)$ suitably along the boundary components corresponding to $p_{(2,3)}$, as follows.

\begin{itemize}
\item If $M(v_3)\cong T^2 \times I$, then  $M(v_2)\cup M(v_3) \cong (S^1 \times D^2) \sharp (S^1 \times D^2)$.
\item If $M(v_3)\cong (S^1 \times D^2) \sharp (S^1 \times D^2)$, then
by the genus-$1$ Heegaard splitting of $S^3$, a suitable gluing gives $M(v_2)\cup M(v_3) \cong (S^1 \times D^2) \sharp (S^1 \times D^2)$.


\item  
If $M(v_3)\cong M(0,2;\frac{1}{2})$, then gluing $M(v_3)$ with a solid torus  by matching regular fibers on one boundary component of $M(v_3)$ to longitudes of the solid torus produces a new solid torus. It follows that a suitable gluing gives $M(v_2)\cup M(v_3) \cong (S^1 \times D^2) \sharp (S^1 \times D^2)$.

\end{itemize}

 Then we always have   $M(v_2)\cup M(v_3)\cong (S^1 \times D^2) \sharp (S^1 \times D^2)$.
Since $\beta_1(L)=1$, there is a saddle vertex in $\{v_4,\dots, v_m\}$  adjacent to $L_{v_2}\cup_{p_{(2,3)}}L_{v_3}$. Repeating the same argument, we have $\cup_{k=2}^{m} M(v_k)\cong (S^1 \times D^2) \sharp (S^1 \times D^2)$.
Recall that $M(v_1) \cong (S^1 \times D^2) \sharp (S^1 \times D^2)$. 
Then by gluing  $\cup_{k=2}^{m} M(v_k)$ and $M(v_1)$ along their boundary suitably, we can obtain \(S^1 \times S^2\).
 
 \end{proo}
 
 By Lemma \ref{lem:block} and Claim \ref{c:inseparable torus gluing}, $L $ can be realized as a WLG of an NMS flow on $S^1\times S^2$. According to the construction,   there exists an inseparable regular level surface  homeomorphic to $T^2$.
  Sufficiency is proved.

\end{proo}

\begin{theorem}\label{thm:inseparable not torus}
 Let $L$ be an  AWLG. $L$  is associated with an NMS flow $\phi_t$ on $S^1 \times S^2$ 
 such that there exists an inseparable regular level surface  not homeomorphic to $T^2$
 if and only if the following conditions hold:
 \begin{enumerate}
   \item $\beta_1(L)=1$; 
   \item 
every edge in the cycle $C$ of $L$
has weight $0$ or $2$, and every edge not  in $C$ has
weight $1$;
\item each saddle vertex is labeled with one of the types listed in 
    Proposition \ref{prop:N-classification} and Proposition \ref{fil nghd twist},
    except for types (B.7), (B.8), and (C.3).

 \end{enumerate}

 \end{theorem}

\begin{proo}
 \emph{Necessity}.  Suppose that $L$ is a WLG of an NMS flow $\phi_t$ on $S^1 \times S^2$ such that  there exists an inseparable regular level surface   not homeomorphic to $T^2$.
 Let $v$ be a saddle vertex of $L$, and  $N $ be the filtrating neighbirhood corresponding to the star neighborhood of $v$.
 
 By Lemma  \ref{lem. inseparable}, $\beta_1(L)=1$, and there exist at least one regular level surface homeomorphic to $S^2$, and one regular level surface homeomorphic to $2T^2$.
 By  Proposition \ref{rek}, every edge  in the unique cycle $C$ of $L$
has weight $0$ or $2$, and every edge not   in $C$ has
weight $1$. Thus there exists an inseparable regular level surface  $\Sigma$ homeomorphic to $2T^2$. 
  By Proposition \ref{prop:2T2}, $S^1 \times S^2 -\Sigma$ can  be embedded in $S^3$,  which implies that  $N$ can be embedded in $S^3$. Therefore, $v$ is labeled with one of the types listed in 
    Proposition \ref{prop:N-classification} and Proposition \ref{fil nghd twist},
    except for types (B.7), (B.8), and (C.3).
\\ 
\par \emph{Sufficiency}. 
Let $L$ be an  AWLG satisfying conditions (1)--(3). Since $\beta_1(L)=1$, we let $C$ be the cycle of $L$. Since every edge not  in $C$ has
weight $1$,  each saddle vertex not in $C$ is labeled with one of  (A.3), (B.2),  (C.1),  (C.4) and  (D.1).
Similar to the proof of sufficiency for Theorem \ref{thm:inseparable torus}, we may assume without loss of generality that $\mathcal{P}_1(L)= \varnothing$.
Cutting $L$ along all points of $\mathcal{P}_2(L)$, we obtain some connected components, each of which contains exactly one saddle vertex. 
Let $L_{v}$ be the component containing  a saddle vertex $v$.

Let \(\mathcal{V}(L)\) denote the set of saddle vertices in \(L\), where for each \(v \in \mathcal{V}(L)\), the two edges connected to \(v\) in the cycle of \(L\) have the same weight.
Let $v_1, \dots, v_m$ be the other saddle vertices of $L$. Namely,   $v_j \notin \mathcal{V}(L)$ for each $j=1,\dots, m$.   By Remark \ref{r:S2,2T2boundary}, we have $m\geq 2$ and $m$ is   even. 
This follows from the fact that   $\beta_{1}(L)=1$ and $L$ contains no oriented cycles.
As previously defined,  the saddle-vertex manifold
$M(v_j)$ associated with $v_j$ is homeomorphic to 
$U-\operatorname{int}D^3$, where $U$ is a genus two handlebody, and $D^3$ is a $3$-ball embedded in $U$. 
By Lemma \ref{lem:block},  there exists an NMS flow $\phi_t^{v_j}$ on $M(v_j)$, such that 
$L_{v_j}$ is a WLG of $\phi_t^{v_j}$, where the ends of $L_{v_j}$ related to the cutting points in $\mathcal{P}_2(L)$ correspond to the boundary of $ M(v_j)$ via the Lyapunov
function. 

Suppose that $v\in \mathcal{V}(L)$ is  adjacent to $v_1$. Let  $p$ be a point in  $\mathcal{P}_2(L)$ adjacent to both $v_1$ and $v$.
Then the saddle-vertex manifold $M(v)\cong S^2 \times I$, or  $M(v)$ is the form  of  
$ U-\operatorname{int}U_0$, where $U_0$ is a genus two handlebody embedded in $U$.
Let $\Sigma$ be a boundary component of $M(v)$. 
In particular, if $M(v)$ is the form  of  
$ U-\operatorname{int}U_0$,  then $\Sigma$ is related to the boundary of $U_0$.
By Lemma \ref{lem:block},  we can obtain an NMS flow $\phi_t^{v}$ on $M(v)$, such that 
$L_{v}$ is a WLG of $\phi_t^{v}$, and  
 the end of $L_{v}$ related to the cutting point $p$  correspond to  $\Sigma$ via the Lyapunov
function.
Then by suitably gluing $M(v)$ and $M(v_1)$ along the boundary components related to $p$, we obtain a manifold homeomorphic to $M(v_1)$.

Replace  $L_v$ with an edge, then we get a new AWLG.   Repeating the same argument for this  new  graph, we  finally obtain a AWLG $L'$ satisfying conditions (1)--(3) and $\mathcal{V}(L')= \varnothing$. 
Obviously, if $L'$   is associated with an NMS flow  on $S^1 \times S^2$, then so is $L$. Now, we may assume without loss of generality that  $\mathcal{V}(L)= \varnothing$.

\begin{figure}[htbp]
\centering
\includegraphics[scale=0.75]{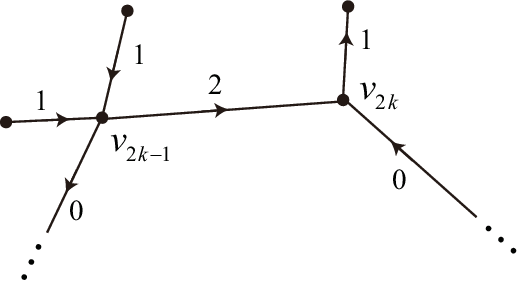}
\caption{ }
\label{Fg_S22T2}
\end{figure}

 Suppose that  $v_{2k-1}$ and $v_{2k}$ are adjacent (see Figure \ref{Fg_S22T2}), and that the adjacent edge has weight $2$ 
 for each $k=1, \dots, \frac{m}{2}$.
 If $m=2$, then by the genus-$2$ Heegaard splitting of $S^3$, we obtain a manifold homeomorphic to $S^2 \times I$ by suitably
  gluing $M(v_1)$ and $M(v_2)$ along their $2T^2$ boundary components.   Gluing the two remaining spherical boundary components together gives
$S^1\times S^2$.
If  $m>2$, then for each $k=1, \dots, \frac{m}{2}$,
 we obtain a manifold homeomorphic to $S^2 \times I$ by suitably
  gluing $M(v_{2k-1})$ and $M(v_{2k})$ along their $2T^2$ boundary components.  
Since $\beta_1(L)=1$, gluing the manifolds  $\{M(v_{2k-1})\cup M(v_{2k})|  k=1, \dots, \frac{m}{2}\}$ according to the gluing rules of $\{L_{v_{2k-1}}\cup L_{v_{2k}}|  k=1, \dots, \frac{m}{2}\}$ in $L$ yields $S^1\times S^2$.

 Therefore,  $L $ can be realized as a WLG of an NMS flow on $S^1\times S^2$. According to the construction,   there exists an inseparable regular level surface not homeomorphic to $T^2$.
  Sufficiency is proved.
\end{proo}

\begin{rema}
Let $L$ be a WLG of  an NMS flow $\phi_t$ on $S^1 \times S^2$. By Lemma \ref{lem. inseparable}, there exists an inseparable regular level surface   not homeomorphic to $T^2$ if and only if there exists a regular level surface  not homeomorphic to $T^2$. Therefore,  Theorem \ref{thm:inseparable not torus} also describes  the Lyapunov graph for which  a  regular level surface not homeomorphic to $T^2$ exists.
\end{rema}

\vskip 1cm
\noindent  Fangfang Chen

\noindent{\small School of Mathematics and Statistics}

\noindent{\small Huangshan University, Huangshan 245041, CHINA}

\noindent{\footnotesize{E-mail: fangfangchen\_97@163.com}}
\vskip 2mm

\end{document}